%% file: manuscript_convergence_PM_arXiv.tex
\documentclass[reqno]{amsart}
\usepackage[top=25truemm,bottom=25truemm,left=25truemm,right=25truemm]{geometry}
\usepackage[dvipdfmx]{graphicx}
\usepackage{amsmath,mathrsfs,amssymb}
\usepackage{enumerate}
\usepackage{here}
\newcommand{\argmin}{\mathop{\rm arg~min}\limits}
\usepackage{cases}
\usepackage{color}
\include{new_command}

\title[Convergence study and optimal weight functions of an explicit particle method]{Convergence study and optimal weight functions of an explicit particle method for the incompressible Navier--Stokes equations}

\author[Y. Imoto]{Yusuke Imoto${}^1$}
\address{${}^1$Kyoto University Institute for Advanced Study, Kyoto University, Yoshida Ushinomiya-cho, Sakyo-ku, Kyoto 6068501, Japan}
\email{imoto.yusuke.4e@kyoto-u.ac.jp}
\author[S. Tsuzuki]{Satori Tsuzuki${}^{2}$}
\address{${}^{2}$Research Center for Advanced Science and Technology, The University of Tokyo, Tokyo, Japan}
\author[D. Nishiura]{Daisuke Nishiura${}^{3}$}
\address{${}^{3}$Department of Mathematical Science and Advanced Technology, Japan Agency For Marine-Earth Science and Technology, Yokohama, Japan}


\begin{document}
	
\begin{abstract}
	To increase the reliability of simulations by particle methods for incompressible viscous flow problems, convergence studies and improvements of accuracy are considered for a fully explicit particle method for incompressible Navier--Stokes equations. 
	The explicit particle method is based on a penalty problem, which converges theoretically to the incompressible Navier--Stokes equations, and is discretized in space by generalized approximate operators defined as a wider class of approximate operators than those of the smoothed particle hydrodynamics (SPH) and moving particle semi-implicit (MPS) methods. 
	By considering an analytical derivation of the explicit particle method and truncation error estimates of the generalized approximate operators, sufficient conditions of convergence are conjectured.
	Under these conditions, the convergence of the explicit particle method is confirmed by numerically comparing errors between exact and approximate solutions. 
	Moreover, by focusing on the truncation errors of the generalized approximate operators, an optimal weight function is derived by reducing the truncation errors over general particle distributions. 
	The effectiveness of the generalized approximate operators with the optimal weight functions is confirmed using numerical results of truncation errors and driven cavity flow. 
	As an application for flow problems with free surface effects, the explicit particle method is applied to a dam break flow. 
\end{abstract}
	
	\maketitle
	
\section{Introduction}
Particle methods, such as the smoothed particle hydrodynamics (SPH) \cite{gingold1977smoothed,lucy1977numerical,price2012smoothed} and moving particle semi-implicit (MPS)  \cite{koshizuka1996moving,khayyer2011enhancement,shakibaeinia2012mps} methods, discretize partial differential equations based on particles distributed in domains and basis functions referred to as weight functions corresponding to each particle. 
These particle methods do not require mesh generation; therefore, they are appropriate for problems that include large deformations or damages, e.g., collapses \cite{monaghan1991simulation}, brittle solids \cite{benz1995simulations}, and Navier--Stokes equations under free surface effects \cite{koshizuka1996moving,liu2003smoothed,monaghan1994simulating,shao2003incompressible}. 
In particular, explicit particle methods for Navier--Stokes equations have been widely used for large-scale problems, such as tsunami run-up \cite{dominguez2013new,murotani2014development}, because of their simple implementation, which can also be done using parallel computing. 

Representative examples of explicit particle methods for the incompressible Navier--Stokes equations include the weekly compressible SPH (WCSPH) \cite{monaghan1994simulating,morris1997modeling} and the explicit MPS (E-MPS) \cite{oochi2010explicit, shakibaeinia2010weakly} methods. 
WCSPH is characterized as an explicit particle method that uses approximate differential operators of SPH for spatial discretization and evaluating pressure given an equation of state for compressible flow. 
In contrast, E-MPS is characterized as an explicit particle method that uses approximate differential operators of MPS for spatial discretization and evaluating pressure given the local density of a number of particles. 
In previous studies \cite{monaghan1994simulating, morris1997modeling, oochi2010explicit, shakibaeinia2010weakly}, both these methods have been validated by comparing their numerical results with experimental results. 
However, in order to ensure reliability for problems on a non-experimental scale or for large-scale computations such as those involving tsunamis, numerical analyses of particle methods, such as convergence studies, are indispensable. 
Although there are a few mathematical analyses for particle methods or related methods \cite{moussa2006convergence,moussa2000convergence,imoto2016phd,imoto2018unique,imoto2019unique,raviart1985analysis}, their results do not directly apply to explicit particle methods. 
Therefore, we focus on convergence studies for explicit particle methods for the incompressible Navier--Stokes equations. 

To describe a mathematical convergence for the incompressible Navier--Stokes equations, in this study, we configure an explicit particle method without physical parameters and assumptions in a manner similar to previous literature  \cite{monaghan1994simulating, morris1997modeling, shakibaeinia2010weakly}. 
Then, we introduce a penalty problem that theoretically converges to the incompressible Navier--Stokes equations and derive the explicit particle method by discretizing the penalty problem based on mathematical theory alone. 
In this spatial discretization, we use generalized approximate operators, which are defined as a wider class of approximate operators for particle methods of SPH and MPS. 
Because of this discretization, computational procedures of the explicit particle methods closely resemble that of E-MPS. 
Furthermore, for the explicit particle method, we conjecture sufficient conditions of convergence based on its analytical derivation and truncation error estimates for the generalized approximate operators. 
Under these sufficient conditions, we confirm the convergence of the explicit particle method by computing errors between numerical solutions and exact solutions in the Taylor--Green vortex.  

Moreover, to improve the accuracy of the explicit particle method, we consider an optimization of the discrete parameters based on the truncation error estimates of the generalized approximate operators \cite{imoto2016phd,imoto2016teintpV,imoto2017tederiV}. 
In particular, defining the generalized approximate operators as a wider class of those used in particle methods enables us to consider an optimization of discrete parameters without imposed constraint conditions in each method. 
Thus, using truncation errors based on particle distributions as the objective function, we introduce an optimization problem for weight functions of the generalized approximate operators. 
The effects of weight functions obtained as solutions of the optimization problem are confirmed by numerical results of truncation errors and a driven cavity flow. 

Furthermore, to confirm that the explicit particle method can be applied to more realistic problems, we develop it for flow problems under free surface effects. 
In the case of the original procedure of the explicit particle method, pressure around a free surface are evaluated as much lower than that in the inner domain of the fluid, because of the lack of particles. 
In addition, clustering of particles around free surface using pressure gradients causes unstable motion. 
Therefore, by modifying the procedure of evaluating pressure and its gradient, we ensure stable simulations of flow problems under free surface effects. 
Moreover, we apply the explicit particle method with these modifications to a dam break flow and compare the obtained numerical and experimental results.

\section{Explicit particle method for incompressible Navier--Stokes equations}
\label{sec:justification_explicit_scheme}
In this section, we present the formulation of the governing equations and approximate operators, which are used for spatial discretization in our study; furthermore, we introduce an explicit particle method for incompressible Navier--Stokes equations. 

\subsection{Governing equation}
\label{subsec:Governing_equation}
Let $\dR$ be the set of real numbers. 
Let $\Dm$ be a bounded domain in $\dRd\,(\dim=2, 3)$ with a smooth boundary $\bd$. 
We consider the incompressible Navier--Stokes equations as follows: 
\begin{subequations}
	\begin{numcases}
	{}
	\mderiv{1}{t}{\Vel} =-\frac{1}{\rho}\grad\Pres+ \nu\lap\Vel + \force,\quad &$(x,t)\in\Dm\times(0,\tmax)$, 
	\label{eq:N-S:moment}\\
	\diver\Vel = 0, \quad&$(x,t)\in\Dm\times(0,\tmax)$, 
	\label{eq:N-S:cont}\\
	\Vel= \Velini, \quad &$x\in\Dm,~t=0$, 
	\label{eq:N-S:ini}\\
	\Vel = \Velbd, \quad&$(x,t)\in\bd\times (0,\tmax)$, 
	\label{eq:N-S:bd}
	\end{numcases}
	\label{eq:N-S}
\end{subequations}
where $\Vel:\Dm\times (0,\tmax)\ra\dRd$, $\Pres:\Dm\times(0,\tmax)\ra\dR$, $\rho>0$, $\nu>0$, $f:\Dm\times(0,\tmax)\ra\dRd$, $\Vel_0:\Dm\ra\dRd$, and $\Velbd:\Gamma\times(0,\tmax)\ra\dRd$ denote velocity, pressure, density, kinematic viscosity, body force, initial velocity, and boundary velocity of the fluid, respectively. 
Furthermore, $\mderivD/\mderivD t$ denotes the material derivative defined as $\mderivD/\mderivD t\deq\partial/\partial t+\Vel\cdot\nabla$. 
The unknown values are velocity $\Vel$ and pressure $\Pres$. 
We assume the uniqueness and existence of a smooth solution for the incompressible Navier--Stokes equations \eqref{eq:N-S}. 
Note that we only treat the Dirichlet boundary condition in \eqref{eq:N-S} for simplicity, although we will deal with boundaries including the free surface in an applied example in Section 5. 

\subsection{Generalized approximate operators}
\label{subsec:Generalized_approximate_operator}
We introduce approximate operators for spatial discretization of an explicit particle method. 
To ensure generality, we use different formulations from those used in the specific cases of SPH and MPS. 

\begin{figure}[t]
	\centering
	\includegraphics[width=55mm,bb=13 17 1837 1185]{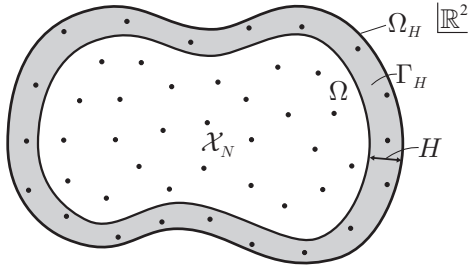}
	\caption{Example of the particle distribution.}
	\label{fig:particle_distribution}
\end{figure}

For a fixed positive number $\H$ and domain $\Dm\subset\dRd$, an expanded domain $\DmH$ is defined as
\begin{align}
	\DmH := \Set{x\in\dRd}{\exists y\in\Dm \mbox{~s.t.~} |x-y|<\H}. 
\end{align}
Let $\bdh\deq\DmH\setminus\Dm$. 
Let $\dN$ be the set of positive integers. 
For $\N\in\dN$, we define a particle distribution $\PtSet$ and particle volume set $\PvSet$ as
\begin{align}
\PtSet&\deq\Set{\Pt{\i}{}\in\DmH}{\i=1,2,\dots,\N,\,\Pt{\i}{}\neq\Pt{\j}{}\,(\i\neq\j)}, 
\label{def:PtSet}\\
\PvSet&\deq\Set{\Pv{\i}>0}{\i=1,2,\dots,\N,\,\sum_{\i=1}^\N\Pv{\i}=\meas{\DmH}}, 
\label{def:PvSet}
\end{align}
respectively. 
Here, $\meas{\DmH}$ indicates the volume of $\DmH$. 
We refer to $\Pt{\i}{}\in\PtSet$ and $\Pv{\i}\in\PvSet$ as a particle and particle volume, respectively.   
Figure \ref{fig:particle_distribution} shows an example of a particle distribution. 

We define a function space $\FsWeightFunc$ as 
\begin{equation}
\FsWeightFunc\deq\Set{\w:[0,\infty)\ra\dR}{\w(r)>0\,(0<r<1),~\w(r)=0\,(r\geq 1)}. 
\label{def:functionspace_w}
\end{equation}
We refer to $\w\in\FsWeightFunc$ as a reference weight function. 
The influence radius $\h$ is a real number satisfying $\min\SetNd{|\Pt{\i}{}-\Pt{\j}{}|}{\i\neq\j}<\h<\H$. 
For the reference weight function $\w$ and influence radius $\h$, a weight function $\wh:[0,\infty)\ra\dR$ is defined as 
\begin{equation}
\wh(r) \deq \dfrac{1}{\h^{\dim}}\w\brs{\dfrac{r}{\h}},\qquad r\in[0,\infty).
\label{def:wh}
\end{equation}
We refer to the domain $\{y\in\dRd;\,|y-\Pt{\i}{}|<h\}$ as an influence domain for particle $\Pt{\i}{}$; in addition, we refer to particles in the influence domain for particle $\Pt{\i}{}$ as the neighbor particles of particle $\Pt{\i}{}$. 
For an integer $k$ and function $\w:[0,\infty)\ra\dR$, we define $C_k(\w)$ as 
\begin{equation}
C_k(\w) \deq \int_{\dRd} |x|^k \w(|x|)\dx. 
\end{equation}

Set discrete parameters $(\PtSet, \PvSet, \h)$ and reference weight functions $\wi, \wg, \wl\in\FsWeightFunc$. 
Then, for $\f:\PtSet\ra\dR$, we define the interpolant $\Intp{}$, approximate gradient operators $\GradApp{}$, and approximate Laplace operator $\LapApp{}$ as
\begin{align}
\Intp{}\f_\i&\deq\IntpCof\sum_{\j=1}^\N\Pv{\j}\f_\j\wih (|\Pt{\j}{}-\Pt{\i}{}|),
\label{approx_intp_01}\\
\GradApp{} \f_\i &\deq \frac{\GradAppCof}{\h} \sum_{\j\neq\i}\Pv{\j} (\f_\j-\f_\i)\frac{\Pt{\j}{}-\Pt{\i}{}}{|\Pt{\j}{}-\Pt{\i}{}|} \wgh (|\Pt{\j}{}-\Pt{\i}{}|),
\label{approx_grad_01}\\
\LapApp{} \f_\i &\deq \frac{\LapAppCof}{\h^2} \sum_{\j\neq\i}\Pv{\j}(\f_\j-\f_\i) \wlh (|\Pt{\j}{}-\Pt{\i}{}|), 
\label{approx_lap_01}
\end{align}
respectively. 
Here, $\f_\i\deq\f(\Pt{\i}{})$, $\IntpCof\deq1/C_0(\wi)$, $\GradAppCof\deq \dim/C_1(\wg)$, and $\LapAppCof\deq 2\dim/C_2(\wl)$. 

The derivations of these operators are presented in Section \ref{subsec:derivation_operator}. 
Moreover, as discussed later in Appendix \ref{sec:appendix_Derivation_conventional}, these operators represent a wider class of approximate operators for particle methods that those in the SPH and MPS methods. 
Thus, we refer to these operators as generalized approximate operators. 
Note that different symbols in reference weight functions $\wi, \wg, \wl$ for each differential operator are used in order to allow us to choose them arbitrarily. 
In Section \ref{subsec:optimal_weight_func}, we discuss the optimization of weight functions for these generalized approximate operators based on truncation error estimates. 

\subsection{Computational procedure of the explicit particle method}
\label{subsec:explicit_scheme}
We introduce an explicit particle method for the incompressible Navier--Stokes equations. 
Before introducing this method, we introduce some notations used in our study. 
Let $\Veliniex:\DmH\ra\dRd$ and $\Velbdex:\bdh\times[0,T]\ra\dRd$ be expanded functions of the initial and boundary velocities, respectively, that satisfy $\Veliniex|_{\Dm} = \Velini$, $\Velbdex|_{\bd} = \Velbd$, and $\Velbdex|_{t=0}=\Veliniex|_{\bdh}$. 
In addition, let $\Dt>0$ be the time step. 
Further, let $\tstepmax$ be the total number of time steps defined by $\tstepmax\deq\gausssymbol{\tmax/\Dt}$, where $\gausssymbol{a}$ denotes the greatest integer that is less than or equal to $a$; this symbol is known as the Gauss symbol. 
For $\tstep=0,1,\dots,\tstepmax$, the $k$th time $\tk$ is defined as $\tk\deq\tstep\,\Dt$. 
Let $\PtSetTime{\k}$ and $\Pt{\i}{\k}$ be a particle distribution and an $\i$th particle in that distribution at $\tk$, respectively. 
Let $\PtSetTent{\k}$ and $\PtP{\i}{\k}$ be a tentative particle distribution and a tentative $\i$th particle in that distribution at $\tk$, respectively. 
For $\w\in\FsWeightFunc$, we define $C_{0,\h}(\w)$, which is an approximation of $C_{0}(\w)$, as 
\begin{equation}
C_{0,\h}(\w) \deq \frac{|\DmH|}{\N} \sum_{z\in\mathbb{Z}^{\dim}} \wh(|\DmH|^{1/\dim}\N^{-1/\dim}\BrA{z}), 
\end{equation}
where $\dZ$ is the set of integers. 
For $S\subset\dRd$, let $\IndexSet{S}{\k}$ be an index set of particles in $S$:
\begin{equation}
\IndexSet{S}{\k}\deq\left\{\i=1,2,\dots,\N\midd\Pt{\i}{\k}\in S\right\}. 
\end{equation}
We denote $\Intp{}$, $\GradApp{}$, and $\LapApp{}$ by replacing $\Pt{\i}{}\in\PtSet$ with $\Pt{\i}{\k}\in\PtSetTime{\k}$ as $\Intp{\k}$, $\GradApp{\k}$, and $\LapApp{\k}$, respectively.  
For $\f:\PtSet\ra\dR$, we define a modified interpolant $\IntpMod{\k}$ and an additional approximate gradient operator $\GradAppPlus{\k}$ as
\begin{align}
\IntpMod{\k}\f_\i&\deq\cfrac{\ds \sum_{\j=1}^{\N} \Pv{\j}\f_j \wih(|\Pt{\j}{\k}-\Pt{\i}{\k}|)}{\ds \sum_{\j=1}^{\N}\Pv{\j}  \wih(|\Pt{\j}{\k}-\Pt{\i}{\k}|)},
\label{approx_intp_02}\\
\GradAppPlus{\k} \f_\i&\deq\frac{\GradAppCof}{\h} \sum_{\j\neq\i}\Pv{\j} (\f_\j+\f_\i)\frac{\Pt{\j}{\k}-\Pt{\i}{\k}}{|\Pt{\j}{\k}-\Pt{\i}{\k}|} \wgh (|\Pt{\j}{\k}-\Pt{\i}{\k}|),
\label{approx_grad_02}
\end{align}
respectively. 
Note that the modified interpolant $\IntpMod{\k}$ corresponds to an interpolant used for a Shepard filter \cite{panizzo2004physical,rogers2010simulation}.

\begin{figure}[tb]
	\centering
	\includegraphics[bb=0 0 186mm 266mm, width=80mm]{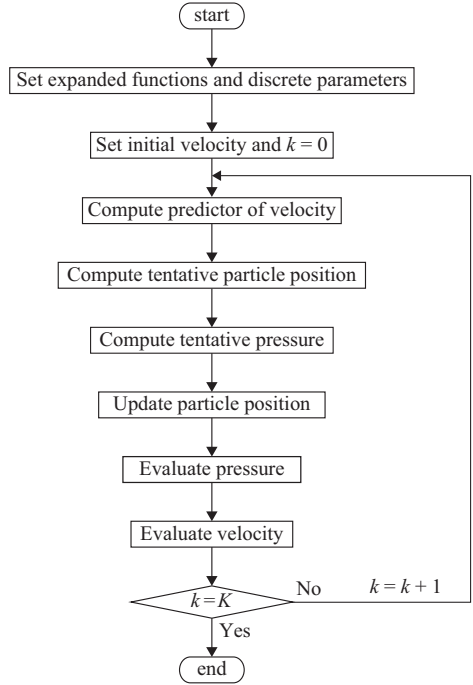}
	\caption{Flowchart of the explicit particle method.}
	\label{fig_flowchart}
\end{figure}

The computational procedure of the explicit particle method involves the following steps. 
Set the discrete parameters as follows: $H>0$, initial particle distribution $\PtSetTime{0}$, particle volume set $\PvSet$, reference weight function $\w$, influence radius $\h\leq\H/2$, parameter $\CompCof>0$, and time step $\Dt$. 
Set the initial approximate velocity $\VelApp{\i}{0}\,(\i=1,2,\dots,\N)$ as $\VelApp{\i}{0}=\Veliniex(\ptt{\i}{0})$. 
Then, for $\tstep=0,1,\dots,\tstepmax-1$, the approximate solution $(\VelApp{\i}{\k+1},\PresApp{\i}{\k+1})\,(\i=1,2,\dots,\N)$ is solved using the following steps:\\
Step 1: Compute a predictor of velocity $\VelAppP{\i}{\k+1}$ as follows:
\begin{equation}
\left\{
\begin{array}{@{\,}r@{\,}c@{\,}ll}
\dfrac{\VelAppP{\i}{\k+1}-\VelApp{\i}{\k}}{\Dt} &=& \nu\LapApp{\k}\VelApp{\i}{\k} +\force(\Pt{\i}{\k},\tk),\qquad&\i\in\IndexSet{\Dm}{\k},
\\
\VelAppP{\i}{\k+1}&=&\Velbdex(\Pt{\i}{\k},\tk),\qquad&\i\in\IndexSet{\bdh}{\k};
\end{array}
\right.
\label{eq:disc_vel_temp}
\end{equation}
Step 2: Compute a tentative particle position $\PtP{\i}{\k+1}$ as follows:
\begin{equation}
\left\{
\begin{array}{@{\,}r@{\,}c@{\,}ll}
\PtP{\i}{\k+1} &=& \Pt{\i}{\k}+\Dt\,\VelAppP{\i}{\k+1},\qquad&\i\in\IndexSet{\Dm}{\k},
\\
\PtP{\i}{\k+1} &=& \Pt{\i}{\k},\qquad&\i\in\IndexSet{\bdh}{\k};
\end{array}
\right.
\label{eq:disc_pd_temp}
\end{equation}
Step 3: Compute a tentative pressure $\PresAppP{\i}{\k+1}$ as follows:
\begin{equation}
\PresAppP{\i}{\k+1} =\dfrac{\rho}{\CompCof^{2}}\brs{\dfrac{1}{C_{0,h}(\w)}\sum_{\j=1}^{\N}\Pv{\j}\wh(|\PtP{\j}{\k+1}-\PtP{\i}{\k+1}|)-1},\qquad\i\in\IndexSet{\DmH}{\k};
\label{eq:disc_pres_temp}
\end{equation}
Step 4: Update the particle position $\Pt{\i}{\k+1}$ as follows:
\begin{equation}
\left\{
\begin{array}{@{\,}r@{\,}c@{\,}ll}
\Pt{\i}{\k+1} &=& \ds \PtP{\i}{\k+1} -\frac{\Dt^2}{\Dens} \GradAppPlustent{\k+1}\PresAppP{\i}{\k+1},\quad&\i\in\IndexSet{\Dm}{\k},
\\
\Pt{\i}{\k+1} &=& \PtP{\i}{\k+1},\qquad&\i\in\IndexSet{\bdh}{\k},
\end{array}
\right.
\label{eq:disc_pd}
\end{equation}
where $\GradAppPlustent{\k+1}$ is the gradient operator $\GradAppPlus{\k+1}$ wherein $\{\Pt{\i}{\k+1}\}$ is replaced with $\{\PtP{\i}{\k+1}\}$; \\
Step 5: Evaluate the pressure $\PresApp{\i}{\k+1}$ as follows:
\begin{equation}
\PresApp{\i}{\k+1} = \IntpMod{\k+1}\PresAppP{\i}{\k+1},\qquad\i\in\IndexSet{\DmH}{\k}.
\label{eq:disc_pres}
\end{equation}
Step 6: Evaluate the velocity $\VelApp{\i}{\k+1}$ as follows:
\begin{equation}
\left\{
\begin{array}{@{\,}r@{\,}c@{\,}ll}
\dfrac{\VelApp{\i}{\k+1}-\VelAppP{\i}{\k+1}}{\Dt} &=& -\dfrac{1}{\rho}\GradApp{\k+1}\PresApp{\i}{\k+1},\qquad&\i\in\IndexSet{\Dm}{\k},
\\
\VelApp{\i}{\k+1}&=&\Velbdex(\Pt{\i}{\k+1},\tkp),\qquad&\i\in\IndexSet{\bdh}{\k}.
\end{array}
\right.
\label{eq:disc_vel}
\end{equation}
The flowchart of the explicit particle method is shown in Figure \ref{fig_flowchart}. 

Because the pressure $\PresApp{\i}{\k+1}$ is evaluated based on the density of neighbor particles, the explicit particle method is similar to E-MPS \cite{shakibaeinia2010weakly}. 
In the next section, we derive the explicit particle method and consider its errors; in addition, we show the convergence of the explicit particle method numerically.

\section{Convergence study}
\label{sec:Derivation_Convergence}
In order to confirm the convergence of the explicit particle method, we conjecture conditions of convergence by considering the truncation error estimates of generalized approximate operators and the derivation of the explicit particle method. 
Moreover, we show the convergence of the explicit particle method using numerical results. 
Note: See Appendix \ref{sec:appendix_notation} for computational rules of the multi-index and definitions of functional spaces and their norms. 

\subsection{Derivation of generalized approximate operators}
\label{subsec:derivation_operator}
In order to estimate truncation errors, we present the derivations of the generalized approximate operators in Section \ref{subsec:Generalized_approximate_operator}. 
Let $\Pt{\i}{}\in\PtSet\cap\Dm$. 
Let $B_r(x)$ be an open ball with the center at $x$ and radius $r$:
\begin{align}
B_r(x)\deq\{y\in\dRd;~|x-y|<r\}. 
\end{align}
Then, by Taylor expansion, for $y\in B_h(\Pt{\i}{})\setminus\{\Pt{\i}{}\}$ and $v\in C^n(\DmHo)~(n\in\dN)$, we have
\begin{align}
\f(y) = \sum_{0\leq|\alpha|\leq n-1} \frac{D^{\alpha}\f_\i}{\alpha!}(y-\Pt{\i}{})^{\alpha} + R_{n,\i}(y;\f). 
\label{taylor01}
\end{align}
Here, $\alpha$ is a multi-index and $R_{n,\i}(y;\f)$ is the residual given by
\begin{align}
R_{n,\i}(y;\f) \deq & \sum_{|\alpha|=n}(y-\Pt{\i}{})^{\alpha} \frac{|\alpha|}{\alpha !} \int_0^1 (1-s)^{|\alpha|-1} D^{\alpha}v(s y+(1-s) \Pt{\i}{})\,{\rm d}s.
\label{taylor_residual}
\end{align}
For $\k=1,2,\dots,\dim$ and nonnegative integer $l$, let $\beta_{\k,l}$ be a multi-index such that the $k$th element is $l$, while the others are $0$. 
For $n=2, 3, 4$ and $k=1,2,\dots,\dim$, multiplying both the sides of \eqref{taylor_residual} by
\begin{align}
\frac{\dim(n-2)!}{\h^{n-2}C_{n-2}(\w)}\frac{(y-\Pt{\i}{})^{\,\beta_{\k,n-2}}}{|y-\Pt{\i}{}|^{n-2}}\wh(|y-\Pt{\i}{}|)
\end{align}
and integrating it over $\DmH$, we get
\begin{multline}
\frac{\dim(n-2)!}{\h^{n-2}C_{n-2}(\w)}\int_{\DmH} \f(y)\frac{(y-\Pt{\i}{})^{\,\beta_{\k,n-2}}}{|y-\Pt{\i}{}|^{n-2}}\wh(|y-\Pt{\i}{}|)\dy\\
=\frac{\dim(n-2)!}{\h^{n-2}C_{n-2}(\w)}\sum_{0\leq|\alpha|\leq n-1} \frac{D^{\alpha}\f_\i}{\alpha!}\int_{\DmH}\frac{(y-\Pt{\i}{})^{\alpha+\beta_{\k,n-2}}}{|y-\Pt{\i}{}|^{n-2}}\wh(|y-\Pt{\i}{}|)\dy + E_{\i,\k,n}. 
\label{deriv:gradlap_01}
\end{multline}
Here, $E_{\i,\k,n}$ is 
\begin{equation}
E_{\i,\k,n}\deq \frac{\dim(n-2)!}{\h^{n-2}C_{n-2}(\w)}\int_{\DmH}R_{n,\i}(y;\f)\frac{(y-\Pt{\i}{})^{\,\beta_{\k,n-2}}}{|y-\Pt{\i}{}|^{n-2}}\wh(|y-\Pt{\i}{}|)\dy=\order(\h^2). 
\label{def:E_i-k-n}
\end{equation}
By considering that
\begin{multline}
\int_{\DmH}\frac{(y-\Pt{\i}{})^{\alpha+\beta_{\k,n-2}}}{|y-\Pt{\i}{}|^{n-2}}\wh(|y-\Pt{\i}{}|)\dy \\
= 
\begin{cases}
0,\qquad&\mbox{one or more elements of $\alpha+\beta_{k, n-2}$ are odd},\\
\dfrac{\h^{n-2}C_{n-2}(\w)}{\dim},\qquad&\alpha+\beta_{k, n-2}=\beta_{k, 2},\\
C_{0}(\w),\qquad&n=2, \alpha=0, 
\end{cases}
\label{deriv:gradlap_02}
\end{multline}
for \eqref{deriv:gradlap_01} with $n=2,3$, we obtain 
\begin{equation}
\f_\i=\dfrac{1}{C_0(\w)}\int_{\DmH}\f(y)\wh(|y-\Pt{\i}{}|)\dy+\order(\h^2)\\
\label{Eikn:n=2}
\end{equation}
and
\begin{equation}
\coord{(\grad\f_\i)}{k}=\frac{\dim}{\h C_{1}(\w)}\int_{\DmH}\{\f(y)-\f_\i\}\frac{\coord{(y-\Pt{\i}{})}{k}}{|y-\Pt{\i}{}|}\wh(|y-\Pt{\i}{}|)\dy+\order(\h^2).
\label{Eikn:n=3}
\end{equation}
Moreover, when $n=4$, by
\begin{equation}
\sum_{k=1}^\dim\frac{(y-\Pt{\i}{})^{\,\beta_{\k,2}}}{|y-\Pt{\i}{}|^{2}}=1
\end{equation}
and
\begin{equation}
\sum_{k=1}^\dim\int_{\DmH}(y-\Pt{\i}{})^{\alpha}\wh(|y-\Pt{\i}{}|)\dy=
\begin{cases}
0,\qquad&\mbox{one or more elements of $\alpha$ are odd},\\
\dfrac{\h^{2}C_2(\w)}{\dim},\qquad&\mbox{$|\alpha|=2$ and all elements of $\alpha$ are even}, 
\end{cases}
\end{equation}
we obtain 
\begin{equation}
\lap\f_\i= \frac{2\dim}{\h^2 C_{2}(\w)}\int_{\DmH}\{\f(y)-\f_\i\}\wh(|y-\Pt{\i}{}|)\dy+\order(\h^2). 
\label{Eikn:n=4}
\end{equation}
By \eqref{def:PvSet}, the above integration can be approximated as
\begin{equation}
\int_{\DmH}\f(y)\dy\approx \sum_{\j=1}^\N\Pv{\j}\f_\j. 
\label{approx:integral}
\end{equation}
Therefore, by \eqref{Eikn:n=2} and \eqref{approx:integral}, and replacing $\w$ with $\wi\in\FsWeightFunc$, we derive the generalized interpolant \eqref{approx_intp_01} as follows:
\begin{align}
\f_\i&=\IntpCof\int_{\DmH}\f(y)\wih(|y-\Pt{\i}{}|)\dy+\order(\h^2)
\nn\\
&\approx \IntpCof\sum_{\j=1}^\N\Pv{\j}\f_\j\wih(|\Pt{\j}{}-\Pt{\i}{}|) = \Intp{}\f_\i. 
\end{align}
By \eqref{Eikn:n=3} and \eqref{approx:integral}, and replacing $\w$ with $\wg\in\FsWeightFunc$, we derive the generalized approximate gradient operator \eqref{approx_grad_01} as follows:
\begin{align}
\grad\f_\i&=\frac{\GradAppCof}{\h }\int_{\DmH}\{\f(y)-\f_\i\}\frac{y-\Pt{\i}{}}{|y-\Pt{\i}{}|}\wgh(|y-\Pt{\i}{}|)\dy+\order(\h^2)
\nn\\
&\approx \frac{\GradAppCof}{\h} \sum_{\j\neq\i}\Pv{\j} (\f_\j-\f_\i)\frac{\Pt{\j}{}-\Pt{\i}{}}{|\Pt{\j}{}-\Pt{\i}{}|} \wgh (|\Pt{\j}{}-\Pt{\i}{}|) = \GradApp{} \f_\i. 
\end{align}
Moreover, by \eqref{Eikn:n=4} and \eqref{approx:integral}, and replacing $\w$ with $\wl\in\FsWeightFunc$, we derive the generalized approximate Laplace operator \eqref{approx_lap_01} as follows:
\begin{align}
\lap\f_\i&=\frac{\LapAppCof}{\h^2}\int_{\DmH}\{\f(y)-\f_\i\}\wlh(|y-\Pt{\i}{}|)\dy+\order(\h^2)
\nn\\
&\approx \frac{\LapAppCof}{\h^2} \sum_{\j\neq\i}\Pv{\j} (\f_\j-\f_\i)\wlh(|\Pt{\j}{}-\Pt{\i}{}|) = \LapApp{} \f_\i. 
\label{derive:laph}
\end{align}
The generalized approximate operators can be used as approximate operators of the conventional particle methods such as SPH and MPS by selecting the parameters of the generalized approximate operators appropriately; this is discussed further in Appendix \ref{sec:appendix_Derivation_conventional}. 
Therefore, approximate operators of conventional particle methods can be derived using the abovementioned method. 

\subsection{Truncation errors of generalized approximate operators}
\label{subsec:truncation_error}
We analyze the truncation errors of the generalized approximate operators using their derivations. 
Let us consider a truncation error estimate of the generalized approximate Laplace operator \eqref{approx_lap_01}. 
We assume $\Pt{\i}{}\in\PtSet\cap\Dm$, $\f\in C^4(\DmHo)$, and  $\wl\in\FsWeightFunc\cap C^1([0,\infty))$. 
From the derivation of the generalized approximate Laplace operator \eqref{derive:laph}, we estimate its truncation error as
\begin{equation}
|\lap\f_\i-\LapApp{}\f_\i| \leq |\widetilde{E}_{\i}| + |\widehat{E}_{\i}|. 
\label{truncat_err:proof01}
\end{equation}
Here,
\begin{equation}
\widetilde{E}_{\i} \deq 2\sum_{k=1}^\dim E_{\i,\k,4} = 2\frac{\LapAppCof}{\h^{2}}\int_{\DmH}R_{4,\i}(y;\f)\wlh(|y-\Pt{\i}{}|)\dy=\order(\h^2),
\label{truncat_err:proof02}
\end{equation}
\begin{equation}
\widehat{E}_{\i} \deq \frac{\LapAppCof}{\h^2}\int_{\DmH}\{\f(y)-\f_\i\}\wlh(|y-\Pt{\i}{}|)\dy - \frac{\LapAppCof}{\h^2} \sum_{\j=1}^\N\Pv{\j} (\f_\j-\f_\i)\wlh(|\Pt{\j}{}-\Pt{\i}{}|). 
\label{truncat_err:proof03}
\end{equation}
Note that the estimate $\widetilde{E}_{\i}=\order(\h^2)$ is derived from \eqref{def:E_i-k-n}. 
	Now, we estimate the error $\widehat{E}_{\i}$, which consists of the integration and the numerical integration, which are the first and second terms on the right hand side of \eqref{truncat_err:proof03}, respectively. 
	For a $C^1$ class function $g:\DmH\rightarrow\dR$ and generators $y_i\in\DmH\,(i=1,2,\dots,N)$, we assume a numerical integration for the integration of $g$ over $\DmH$ given by
	\begin{equation}
	\sum_{i=1}^N|\sigma_i|g(y_i). 
	\end{equation}
	Here, $\sigma=\{\sigma_i\}_{\i=1}^{\N}$ is a decomposition of $\DmH$ satisfying
	\begin{equation}
	\bigcup_{\i=1}^{\N} \overline{\sigma}_\i = \DmHo,
	\qquad
	\sigma_i \cap \sigma_j =\emptyset\quad (i\neq j),
	\end{equation}
	where $\overline{\sigma}_\i$ is the closure of $\sigma_\i$. 
	Then, as an estimate of the Riemann sum, we can estimate the numerical integration as
	\begin{equation}
	\left| \int_{\DmH} g(y)dy -\sum_{i=1,2,\dots,N}|\sigma_i|g(y_i)\right| =\order\left( \max_{i=1,2,\dots,N}{\rm rad}(\sigma_i)\right).
	\end{equation}
	Here, ${\rm rad}(\sigma_i):=\sup\Set{|y_i-z|}{z\in\sigma_i}$. 
	Furthermore, because $\sigma$ is arbitrary, 
	we can estimate the numerical integration as
	\begin{equation}
	\left| \int_{\DmH} g(y)dy -\sum_{i=1,2,\dots,N}|\sigma_i|g(y_i)\right| =\order\left(\inf_{\sigma}\max_{i=1,2,\dots,N}{\rm rad}(\sigma_i)\right). 
	\end{equation}

From the strategy above, we introduce a decomposition of $\DmH$ as $\sigma=\{\sigma_i\}_{\i=1}^{\N}$ such that
\begin{align}
|\sigma_\i| = \Pv{\i}\quad(\i=1,2,\dots,\N),
\qquad
\bigcup_{\i=1}^{\N} \overline{\sigma}_\i = \DmHo,
\qquad
\sigma_i \cap \sigma_j =\emptyset\quad (i\neq j). 
\label{cond:sigma}
\end{align}
For $\sigma$, indicator $\delta_{\sigma}$ is defined as
\begin{align}
\delta_{\sigma} \deq \max_{\i=1,2,\dots,\N}\max_{x\in\overline{\sigma}_i}|\Pt{\i}{}-x|
\end{align}
and indicator $\delta_\infty=\delta_\infty(\PtSet,\PvSet)$ is defined as
\begin{align}
\delta_\infty \deq \inf_{\sigma} \delta_{\sigma}. 
\end{align}

Let any $\sigma=\{\sigma_i\}_{\i=1}^{\N}$ such that \eqref{cond:sigma}.  
Furthermore, we assume $\delta_\infty\leq\h$. 
Then, by Taylor's theorem, we can estimate the following:
\begin{align}
|\widehat{E}_{\i}| 
&=\frac{\LapAppCof}{\h^2}\BrA{\int_{\DmH}\{\f(y)-\f_\i\}\wlh(|y-\Pt{\i}{}|)\dy-\sum_{\j=1}^\N\Pv{\j} (\f_\j-\f_\i)\wlh(|\Pt{\j}{}-\Pt{\i}{}|)}
\nn\\
&\leq\frac{\LapAppCof}{\h^2}\BrA{\int_{\DmH}\{\f(y)-\f_\i\}\wlh(|y-\Pt{\i}{}|)\dy-\sum_{\j=1}^\N (\f_\j-\f_\i)\int_{\sigma_\j}\wlh(|y-\Pt{\i}{}|)\dy}
\nn\\
&\quad+\frac{\LapAppCof}{\h^2}\BrA{\sum_{\j=1}^\N (\f_\j-\f_\i)\int_{\sigma_\j}\wlh(|y-\Pt{\i}{}|)\dy-\sum_{\j=1}^\N\Pv{\j} (\f_\j-\f_\i)\wlh(|\Pt{\j}{}-\Pt{\i}{}|)}
\nn\\
&\leq\frac{\LapAppCof}{\h^2}\BrA{\sum_{\j=1}^\N \int_{\sigma_\j}\{\f(y)-\f_\j\}\wlh(|y-\Pt{\i}{}|)\dy}+\frac{\LapAppCof}{\h^2}\BrA{\sum_{\j=1}^\N (\f_\j-\f_\i)\int_{\sigma_\j}\{\wlh(|y-\Pt{\i}{}|)-\wlh(|\Pt{\j}{}-\Pt{\i}{}|)\}\dy}
\nn\\
&\leq\frac{\delta_{\sigma}}{\h^2}|\f|_{C^1(\DmHo)}\LapAppCof C_{1}(\wl)+\frac{\h+\delta_{\sigma}}{\h^2}|\f|_{C^1(\DmHo)}\LapAppCof\sum_{\j\in\{\k ;\,|\Pt{\k}{}-\Pt{\i}{}|<\h+\delta_{\sigma}\}}\int_{\sigma_\j}|\wlh(|y-\Pt{\i}{}|)-\wlh(|\Pt{\j}{}-\Pt{\i}{}|)|\dy
\nn\\
&=\frac{\delta_{\sigma}}{\h^2}|\f|_{C^1(\DmHo)}\LapAppCof C_{1}(\wl)+\brs{1+\frac{\delta_{\sigma}}{\h}}\frac{\delta_{\sigma}}{\h}|\f|_{C^1(\DmHo)}\LapAppCof\int_{\dRd}\BrA{\deriv{1}{r}{}\wlh(|y-\Pt{\i}{}|)}\dy+\order(\delta_{\sigma}^2\h^{-3})
\nn\\
&=\LapAppCof\brs{C_{0}(\wl)+2\int_{\dRd}\BrA{\deriv{1}{r}{}\wl(|y|)}\dy}\frac{\delta_{\sigma}}{\h^2}|\f|_{C^1(\DmHo)}+\order(\delta_{\sigma}^2\h^{-3})
\end{align}
Because $\sigma$ is arbitrary, we obtain
\begin{equation}
|\widehat{E}_{\i}| =\LapAppCof\brs{C_{0}(\wl)+2\int_{\dRd}\BrA{\deriv{1}{r}{}\wl(|y|)}\dy}\frac{\delta_{\infty}}{\h^2}|\f|_{C^1(\DmHo)}+\order(\delta_{\infty}^2\h^{-3}). 
\end{equation}
Hence, $\delta_\infty\leq\h$ yields
\begin{equation}
|\widehat{E}_{\i}| = \order(\delta_{\infty}\h^{-2}). 
\label{error_Lia_01}
\end{equation}
Consequently, by \eqref{truncat_err:proof01} and \eqref{truncat_err:proof02}, and \eqref{error_Lia_01}, we establish
\begin{equation}
|\lap\f_\i-\LapApp{}\f_\i| = \order(\h^2+\delta_{\infty}\h^{-2}). 
\end{equation}

Let $r_{\min}$ be $r_{\min}\deq\min\{|\Pt{\j}{}-\Pt{\i}{}|;\,\i,\j=1,2,\dots,\N, \i\neq \j\}$. 
If $\delta_\infty=\order(r_{\min})$, we refer to the particle distribution and particle volume as regular. 
Based on the definition of $\delta_\infty$, in the absence of extremely unfavorable conditions, such as high density particle distributions or high variance particle volumes, the particle distribution and particle volume become regular. 
The indicator $\delta_\infty$ satisfies $\delta_\infty=\order(\N^{-1/\dim})=\order(r_{\min})$. 
By assuming the regularity of the particle distribution and particle volume, we estimate the truncation error of the generalized approximate Laplace operator as 
\begin{equation}
|\lap\f_\i-\LapApp{}\f_\i| = \order(\h^2+r_{\min}\h^{-2}). 
\label{order:truncation_error}
\end{equation}
A more precise theorem to estimate truncation error has been reported in the literature\cite{imoto2016phd,imoto2016teintpV,imoto2016teintpV}. 

\subsection{Derivation of the explicit particle method}
\label{subsec:derivation}
The explicit particle method is based on the following penalty problem for the incompressible Navier--Stokes equations: 
\begin{subequations}
	\begin{numcases}
	{}
	\MDerivComp{\VelComp} =-\frac{1}{\rho}\grad\PresComp+ \nu\lap\VelComp + \force,\quad &$(x,t)\in\Dm\times(0,\tmax)$, 
	\label{eq:N-S_Comp:moment}\\
	\CompCof^{2}\MDerivComp{\PresComp}+\Dens\diver\VelComp = 0, \quad&$(x,t)\in\Dm\times(0,\tmax)$, 
	\label{eq:N-S_Comp:cont}\\
	\VelComp= \Veliniex, \quad &$x\in\DmH,~t=0$, 
	\label{eq:N-S_Comp:ini}\\
	\VelComp = \Velbdex, \quad&$(x,t)\in\bdh\times(0,\tmax)$, 
	\label{eq:N-S_Comp:bd}\\
	\PresComp = \Pres_0, \quad &$x\in\Dm,~t=0$.  
	\label{eq:N-S_Comp:ini_dens}
	\end{numcases}
	\label{eq:N-S_Comp}
\end{subequations}
Here, $\CompCof$ and $\Pres_0$ are a penalty term in $\dR$ and the initial pressure, respectively. 
Furthermore, $\MDerivComp{}$ denotes a material derivative defined as $\MDerivComp{}\deq\partial/\partial t+\VelComp\cdot\nabla$. 
The unknown values include $\VelComp:\DmH\times[0,\tmax)\ra\dRd$ and $\PresComp:\Dm\times[0,\tmax)\ra\dR$. 
\eqref{eq:N-S_Comp:moment} is the moment equation, which is the same as \eqref{eq:N-S:moment}. 
Further, \eqref{eq:N-S_Comp:cont} is based on the continuity equation for the compressible flow. 
If $\CompCof=0$, we find that \eqref{eq:N-S_Comp:cont} is equivalent to \eqref{eq:N-S:cont}. 
Therefore, the solution $(\VelComp,\PresComp)$ in the penalty problem \eqref{eq:N-S_Comp} coincides with the solution $(\Vel,\Pres)$ in the original incompressible Navier--Stokes equations \eqref{eq:N-S} if $\CompCof=0$ formally. 
In particular, in the cases of two-dimensional spaces and partially- or full-periodic boundary conditions, the convergence of the penalty problem \eqref{eq:N-S_Comp} has orders of velocity and pressure as $\order(\CompCof^{2})$ and $\order(\CompCof)$, respectively, which has been proved in Kreiss et al.\cite{kreiss1991convergence}. 

We arbitrarily set $\tstep=0,1,\dots,\tstepmax-1$. 
Before deriving the discretized schemes, we define a function $\PresCompApp{\k}:\Dm\times[\tk,\tkp)\ra\dR$ as
\begin{equation}
\PresCompApp{\k}(x,t)\deq \dfrac{\rho}{\CompCof^{2}}\brs{\dfrac{1}{C_0(\w)} \int_{\DmH}\wh(|\XComp{\k}(y,t)-\XComp{\k}(x,t)|) \dy-1}, 
\label{def:PresCompApp}
\end{equation}
where $\w\in\FsWeightFunc$ and $\XComp{\k}$ is the solution of the following differential equation:
\begin{equation}
\left\lbrace 
\begin{array}{@{\,}r@{\,}c@{\,}ll}
\MDerivComp{\XComp{\k}}(x,t) &=&\VelComp(\XComp{\k}(x,t),t),\qquad &t\in(\tk,\tkp), \\
\XComp{\k}(x,t) &=&x,\qquad & t=\tk. 
\end{array}\right. 
\label{def:XComp}
\end{equation}
Then, under the assumption that $\brnnd{\VelComp}_{C^1([0,\tmax];C^3(\DmH))}<\infty$, we have
\begin{equation}
\CompCof^2\MDerivComp{\PresCompApp{\k}}(x,t)+\Dens\diver\VelComp(x,t) =\order(\tau\h^{-1}+\h^2),\qquad x\in\Dm,\quad t\in[\tk,\tkp).
\label{eq:N-S_Comp:cont_app}
\end{equation}
The proof for which is presented in Appendix \ref{sec:appendix_proof_discrete_density}. 
By comparing \eqref{eq:N-S_Comp:cont} and \eqref{eq:N-S_Comp:cont_app}, the function $\PresCompApp{\k}$ yields an approximation of the pressure $\PresComp$ at $t\in[\tk,\tkp)$. 

Next, we introduce a time-discretized scheme for the penalty problem \eqref{eq:N-S_Comp}. 
For $\k=0,1,\dots,\tstepmax$, let $(\VelCompDt{\k}, \PresCompDt{\k})$ be a solution of this scheme at $t=\tk$. 
We set $\k=1,2,\dots,\tstepmax-1$ and $x\in\Dm$ in an arbitrary manner. 
For $y\in\Dm$, we introduce $\XCompDt{\k+1}(y)\,(\approx\XComp{\k}(y,\tkp))$ as
\begin{equation}
\XCompDt{\k+1}(y) \deq y + \Dt\VelCompDt{\k+1}(\XCompDt{\k+1}(y)). 
\end{equation}
Because the material derivative is estimated by
\begin{equation}
\MDerivComp{\f}(x,t) = \dfrac{\f(\XComp{\k}(x,\tkp),\tkp)-\f(x,\tk)}{\Dt}+ \order(\Dt),\qquad t\in(\tk,\tkp), 
\label{order:MDerivComp}
\end{equation}
we discretize \eqref{eq:N-S_Comp:moment} as 
\begin{equation}
\dfrac{\VelCompDt{\k+1}(\XCompDt{\k+1}(x))-\VelCompDt{\k}(x)}{\Dt} =-\frac{1}{\rho}\grad\PresCompDt{\k+1}(\XCompDt{\k+1}(x))+ \nu\lap\VelCompDt{\k}(x) + \force(x,\tk). 
\label{disc-time:moment}
\end{equation}
To evaluate the approximate pressure $\PresCompDt{\k+1}$, we introduce a tentative velocity $\VelCompDtP{\k+1}$ and tentative position $\XCompDtP{\k+1}$ for $x\in\Dm$ as 
\begin{equation}
\VelCompDtP{\k+1}(\XCompDt{\k+1}(x)) \deq \VelCompDt{\k}(x) + \Dt\brm{\nu\lap\VelCompDt{\k}(x) + \force(x,\tk)}
\label{def:VelCompDtP}
\end{equation}
and 
\begin{equation}
\XCompDtP{\k+1}(x) \deq x + \Dt\,\VelCompDtP{\k+1}(\XCompDt{\k+1}(x)), 
\label{def:XCompDtP}
\end{equation}
respectively. 
Using these equations and \eqref{def:PresCompApp}, the approximate pressure $\PresCompDt{\k+1}$ is obtained as follows:
\begin{equation}
\PresCompDt{\k+1}(x) = \dfrac{\rho}{\CompCof^{2}}\brs{\dfrac{1}{C_0(\w)} \int_{\DmH}\wh(|\XCompDtP{\k+1}(y,t)-\XCompDtP{\k+1}(x,t)|) \dy-1}. 
\label{def:PresCompDt}
\end{equation}

Then, by discretizing the time-discretized scheme in space, we derive the explicit particle method. 
Let $\i=1,2,\dots,\N$ such that $\Pt{\i}{\k}\in\Dm$. 
First, we discretize \eqref{def:VelCompDtP} and \eqref{def:XCompDtP} as 
\begin{equation}
\VelAppP{\i}{\k+1} = \VelApp{\i}{\k} + \Dt\brm{\nu\LapApp{\k}\VelApp{\i}{\k} + \force(\Pt{\i}{\k},\tk)}
\label{def:VelAppP}
\end{equation}
and 
\begin{equation}
\PtP{\i}{\k+1} = \Pt{\i}{\k}+\Dt\,\VelAppP{\i}{\k+1},
\end{equation}
respectively. 
Using these, we discretize \eqref{def:PresCompDt} as
\begin{equation}
\PresAppP{\i}{\k+1} =\dfrac{\rho}{\CompCof^{2}}\brs{\dfrac{1}{C_{0,h}(\w)}\sum_{\j=1}^{\N}\Pv{\j}\wh(|\PtP{\j}{\k+1}-\PtP{\i}{\k+1}|)-1}. 
\end{equation}
Then, the particle position is updated as follows:
\begin{equation}
\Pt{\i}{\k+1} = \PtP{\i}{\k+1}-\dfrac{\Dt^2}{\Dens}\GradAppPlustent{\k+1}\PresAppP{\i}{\k+1}. 
\label{update:Pt}
\end{equation}
In this case, to avoid non-uniform particle distributions as discussed in Price \cite{price2012smoothed}, we use $\GradAppPlustent{\k+1}$ as the gradient operator in \eqref{update:Pt}. 
Furthermore, to eliminate noise corresponding to the density of particles, we modify the pressure calculation as follows:
\begin{equation}
\PresApp{\i}{\k+1} = \IntpMod{\k+1}\PresAppP{\i}{\k+1}.
\label{derivation:EPM:pressure_recalc} 
\end{equation}
Then, by using pressure $\PresApp{\i}{\k+1}$, we discretize \eqref{disc-time:moment} as
\begin{equation}
\dfrac{\VelApp{\i}{\k+1}-\VelApp{\i}{\k}}{\Dt} = -\dfrac{1}{\rho}\GradApp{\k+1}\PresApp{\i}{\k+1} + \LapApp{\k}\VelApp{\i}{\k} + \force(\Pt{\i}{\k},\tk). 
\label{update:VelApp}
\end{equation}
Finally, using \eqref{def:VelAppP} and \eqref{update:VelApp}, we derive \eqref{eq:disc_vel}. 
The pressure recalculation \eqref{derivation:EPM:pressure_recalc} is essential to obtain stable and accurate results as shown in numerical experiments in Section \ref{subsec:taylor_green}. 

\subsection{Sufficient conditions of convergence}
\label{subsec:conjecture}
We conjecture the sufficient conditions of convergence for the explicit particle method by considering the deviations and truncation error estimates that were calculated in previous sections. 
In particular, as sufficient conditions of convergence, we require $\h\ra0$ and $r_{\min}\h^{-2}\ra0$ from the truncation error estimates \eqref{order:truncation_error} calculated in Section \ref{subsec:truncation_error}. 
In addition, because the convergence orders between the solution of the incompressible Navier--Stokes equation \eqref{eq:N-S} and that of the penalty problem \eqref{eq:N-S_Comp} are $\order(\CompCof)$, we require $\CompCof\ra0$.  
Moreover, we require $\tau\h^{-1}\ra0$ and $\Dt\ra0$ from the order estimates \eqref{eq:N-S_Comp:cont_app} and \eqref{order:MDerivComp} obtained in Section \ref{subsec:derivation}. 
Thus, the summary of the above conditions is as follows: $\h\ra0$, $r_{\min}\h^{-2}\ra0$, $\CompCof\ra0$, $\tau\h^{-1}\ra0$, and $\Dt\ra0$. 
In particular, when the time step $\Dt$ satisfies 
\begin{equation}
\Dt \leq \Dtmax \deq \min\left\{\dfrac{\h\CompCof}{4}, \dfrac{\h^{1/2}}{4\|\force\|_{L^\infty([0,T]; L^\infty(\Dm))}^{1/2}}, \dfrac{\h^2}{8\nu}\right\}, 
\label{cond:dt}
\end{equation}
where $\|\cdot\|_{L^\infty([0,T]; L^\infty(\Dm))}$ denotes the infinity norm in space-time $\|\f\|_{L^\infty([0,T]; L^\infty(\Dm))} \deq \esssup\{|\f(x,t)|; and \,x\in\Dm, t\in(0,\tmax)\}$; 
the conditions for convergence are
\begin{equation}
\CompCof\ra 0,\quad \h\ra0,\quad r_{\min}\h^{-2}\ra 0.
\label{cond:convergence}
\end{equation}
Condition \eqref{cond:dt} is based on the von Neumann stability analysis and corresponds to that suggested in Morris et al. \cite{morris1997modeling} by replacing the sound speed with $\CompCof^{-1}$. 
Under the conjectured sufficient conditions, we confirm the convergence of the explicit particle method numerically; this is shown in the next subsection. 

\subsection{Numerical convergence}
\label{subsec:taylor_green}
We confirm the conjectured sufficient condition of convergence using the numerical results for the Taylor--Green vortex. 
The Taylor--Green vortex is one of the solutions of the two-dimensional incompressible Navier--Stokes equations \eqref{eq:N-S} in the absence of body force ($f\equiv0$). 
Let $\Dm=(0,L)\times(0,L)$. 
The solutions of the Taylor--Green vortex $(\Vel=(\Vel_1,\Vel_2)^{T}, \Pres)$ are given by
\begin{align}
\Vel_1(x) &= -Ue^{-8\pi^2 t/Re}\cos(2\pi\coord{x}{1}/L)\sin(2\pi\coord{x}{2}/L),
\label{sol:TGvortex_velocity_x}\\
\Vel_2(x) &= Ue^{-8\pi^2 t/Re}\sin(2\pi\coord{x}{1}/L)\cos(2\pi\coord{x}{2}/L),
\label{sol:TGvortex_velocity_y}\\
\Pres(x) &= -\dfrac{\rho}{4}e^{-16\pi^2 t/Re}\{\cos(4\pi\coord{x}{1}/L)+\cos(4\pi\coord{x}{2}/L)\} . 
\label{sol:TGvortex_pressure}
\end{align}
Here, $U$ is the velocity scale, and $Re$ is the Reynolds number defined as $Re\deq UL/\nu$. 
Hereafter, we set $T=0.1$, $\rho=1$, $U=1$, $L=1$, and $\nu=10^{-1}$, namely, $Re=10$. 
By comparing the exact solution and a numerical solution of the Taylor--Green vortex, we investigate the validity of the accuracy of the pressure recalculation \eqref{eq:disc_pres} and convergences of the explicit particle method. 
It should be noted that we do not treat a comparison of accuracy for approximate operators here because the Taylor--Green vortex represents an isotopic flow and disturbances in particle distributions rarely appear in the case when the explicit particle method is used. 

\begin{figure}[th]
	\begin{center}
		\includegraphics[bb=0 0 291mm 131mm,width=100mm]{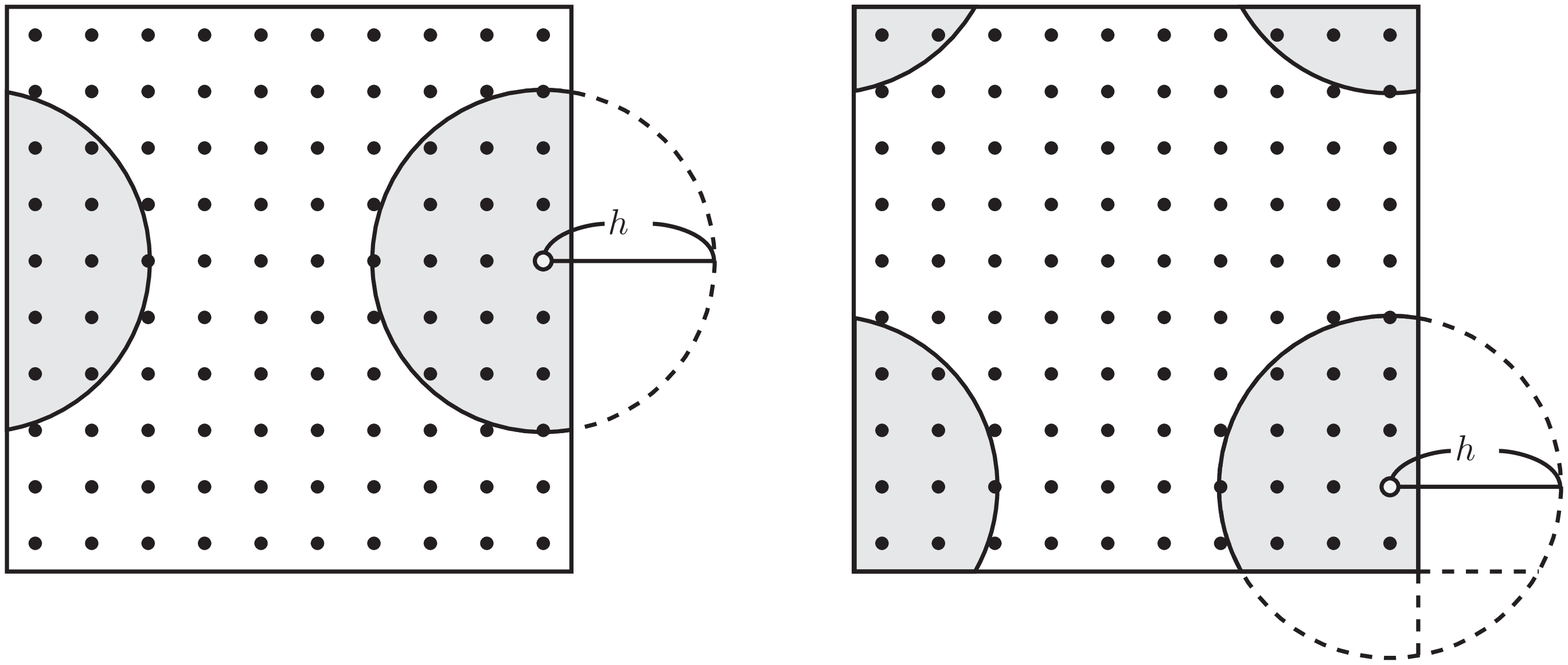}
		\caption{Influence domains near the boundary}
		\label{fig:influence_domain_periodric}
	\end{center}
\end{figure}

\begin{figure}[th]
	\begin{center}
		\includegraphics[bb=0 0 291mm 131mm,width=100mm]{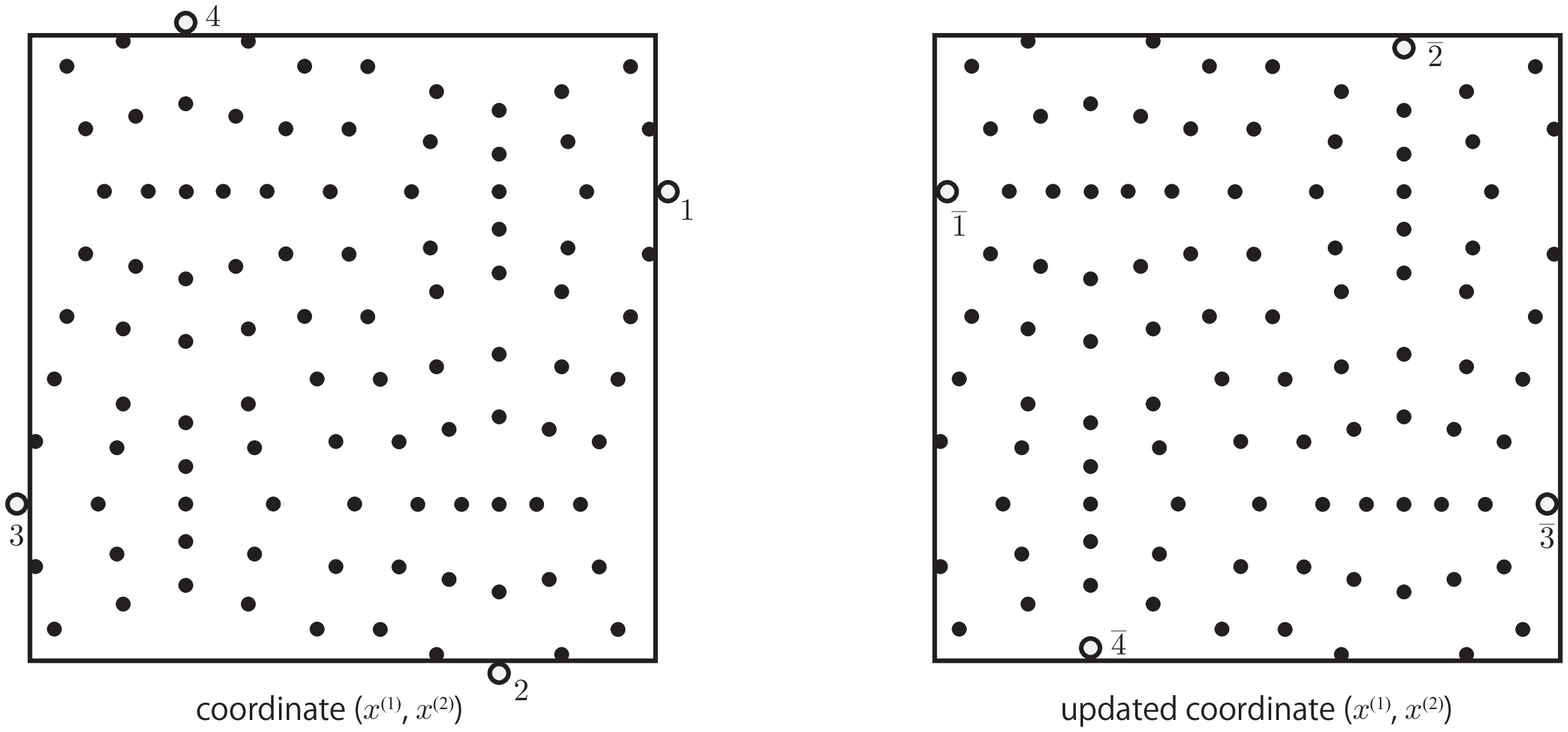}
		\caption{Periodic movements of particles}
		\label{fig:periodric_particle movemnt}
	\end{center}
\end{figure}

Before performing the numerical experiments for convergence, we confirm the computational stability and accuracy of the explicit particle method. 
Because the Taylor--Green vortex is periodic in space, we consider a periodic domain. 
In particular, we consider the following coordinate system: ($\coord{x}{1}, \coord{x}{2}$); $\coord{x}{k}\leftarrow\coord{x}{k}+1$ if $\coord{x}{k}\leq0$ and $\coord{x}{k}\leftarrow\coord{x}{k}-1$ if $\coord{x}{k}\geq1$ for $k=1,2$. 
Then, the particles near the boundary refer to the influence domain corresponding to the periodic boundary conditions as shown in Figure \ref{fig:influence_domain_periodric}.
Moreover, if a particle crosses over the boundary, we let the particle move according to the treatment shown in Figure \ref{fig:periodric_particle movemnt}. 
Because the boundary condition is not required for the system, we do not set the parameter $H$ and expanded boundary condition $\Velbdex$ for it. 
The initial particle distribution $\PtSet^{\timeindex{0}}$ is set as the square lattice with spacing $\Dx=0.04$: 
\begin{equation}
\PtSet^{\timeindex{0}} =\left\{\left(\left(i-\frac{1}{2}\right)\Dx, \left(j-\frac{1}{2}\right)\Dx\right)\in\Dm;\,i,j=1,2,\dots,\lfloor 1/\Dx\rfloor\right\}. 
\label{TGvortex_ini_part_dist}
\end{equation}
Then, the number of particles is $N=(\lfloor 1/\Dx\rfloor)^2=25^2$. 
Furthermore, the particle volume set $\PvSet=\{\Pv{\i}\}_{\i=1}^\N$ is set as
\begin{align}
\Pv{\i} = \dfrac{|\Dm|}{\N}=\Dx^2,\qquad \i=1,2,\dots,\N. 
\label{pt_volume_eg}
\end{align}
We consider five sets of reference weight functions $(\wi, \wg, \wl)$;
\begin{enumerate}
	\setlength{\itemindent}{3.5mm}
	\setlength{\parskip}{2ex}
	\setlength{\itemsep}{0ex}
	\item[(G-s)] $\wi$, $\wg$, and $\wl$ are set as 
		\begin{equation}
		\wi(r)=\w^{\rm spike}(r),\quad
		\wg(r)=\w^{\rm spike}(r), \quad
		\wl(r)=\w^{\rm spike}(r), 
		\label{weight_setting_spike}
		\end{equation}
		where $\w^{\rm spike}$ is the spike function given by
		\begin{align}
		\w^{\rm spike}(r)=
		\begin{cases}
		(1-r)^2,\quad&0\leq r<1,\\
		0,\quad & r\geq 1; 
		\end{cases}
		\label{weight_spike}
		\end{align}
	\item[(S-c)] $\wi$, $\wg$, and $\wl$ are set as 
	\begin{equation}
	\wi(r)=\ws(r),\quad
	\wg(r)=-\WSPHDriv(r), \quad
	\wl(r)=-\dfrac{1}{r}\WSPHDriv(r), 
	\label{weight_SPH}
	\end{equation}
	where $\WSPHDriv$ is the first derivative of $\ws$, in which $\ws$ uses the cubic B-spline defined as
	\begin{equation}
	\w^{\rm cubic}(r) \deq a_{\dim}^{\rm cubic}
	\begin{cases}
	1-6r^2+6r^3, \quad & \displaystyle 0 \leq r < \frac{1}{2},
	\vspace*{0.5ex}\\
	2\left(1-r\right)^3, \quad & \displaystyle \frac{1}{2} \leq r <1,
	\vspace*{0.5ex}\\
	0, \quad & r \geq 1;
	\end{cases}
	\label{w_SPH_cubic}
	\end{equation}
	\item[(S-q)] $\wi$, $\wg$, and $\wl$ are set by \eqref{weight_SPH} in which $\ws$ uses the quintic B-spline defined as
	\begin{equation}
	\w^{\rm quintic}(r) \deq a_{\dim}^{\rm quintic}
	\begin{cases}
	\left(3-3r\right)^5+6\left(2-3r\right)^5+15\left(1-3r\right)^5, \quad & \displaystyle 0 \leq r < \frac{1}{3},
	\vspace*{0.5ex}\\
	\left(3-3r\right)^5+6\left(2-3r\right)^5, \quad & \displaystyle \frac{1}{3} \leq r < \frac{2}{3},
	\vspace*{0.5ex}\\
	\left(3-3r\right)^5, \quad & \displaystyle \frac{2}{3} \leq r <1,
	\vspace*{0.5ex}\\
	0, \quad & r \geq 1;
	\end{cases}
	\label{w_SPH_quintic}
	\end{equation}
	\item[(S-w)] $\wi$, $\wg$, and $\wl$ are set by \eqref{weight_SPH} in which $\ws$ uses the quintic Wendland function (a fifth positive definite function) defined as
	\begin{equation}
	\w^{\rm Wendland}(r) \deq a_{\dim}^{\rm Wendland}
	\begin{cases}
	\left(1-r\right)^4\left(1+4r\right), \quad & 0 \leq r < 1,
	\\
	0, \quad & r \geq 1; 
	\end{cases}
	\label{w_SPH_wendland}
	\end{equation}
	\item[(M)] $\wi$, $\wg$, and $\wl$ are set as
	\begin{equation}
	\wi(r)=\wm(r),\quad
	\wg(r)=\dfrac{1}{r}\wm(r),\quad
	\wl(r)=\wm(r),
	\label{weight_MPS}
	\end{equation}
	respectively, where $\wm$ is the reference weight function of MPS defined as
	\begin{align}
	\wm(r) &\deq
	\begin{cases}
	\dfrac{1}{r}-1, \quad& 0<r<1,\\
	0, \quad& r=0, r\geq 1. 
	\end{cases}
	\label{def:w_MPS}
	\end{align}
\end{enumerate}
Here, $a_{\dim}^{\rm cubic}$, $a_{\dim}^{\rm quintic}$, and $a_{\dim}^{\rm Wendland}$ are constants that satisfy the unity condition:
\begin{equation}
\int_{\dRd}\ws(|x|)dx=1.
\end{equation} 
As shown in Appendix \ref{sec:appendix_Derivation_conventional}, the cases (S-c), (S-q), and (S-w) correspond to the use of approximate operators in SPH. 
Further, case (M) corresponds to the use of approximate operators in MPS.
The influence radius is set as $\h=0.124\,(=3.1\Dx)$. 
In addition, we set $\CompCof=0.1$ and $\Dt=\Dtmax$. 

\begin{figure}[th]
	\flushleft \qquad\quad (a) velocity\\
	\centering
	\includegraphics[bb=0 0 226.3mm 108.0mm, width=130mm]{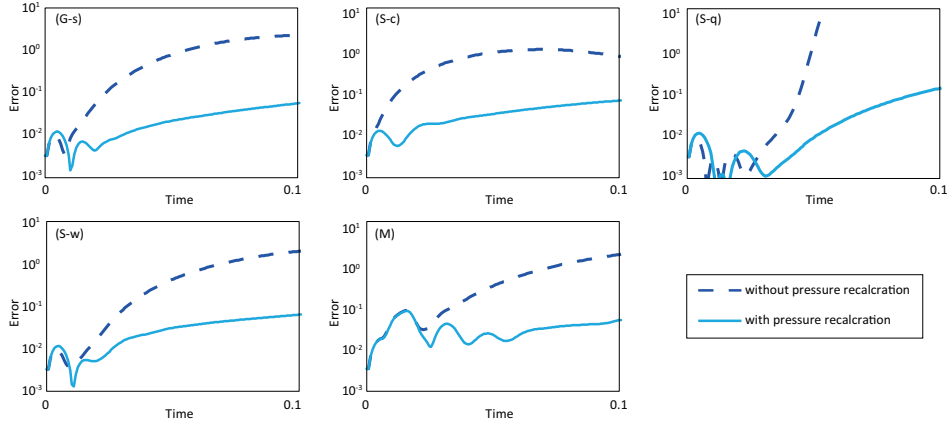}
	\flushleft \qquad\quad (b) pressure\\
	\centering
	\includegraphics[bb=0 0 220.21mm 101.95mm, width=130mm]{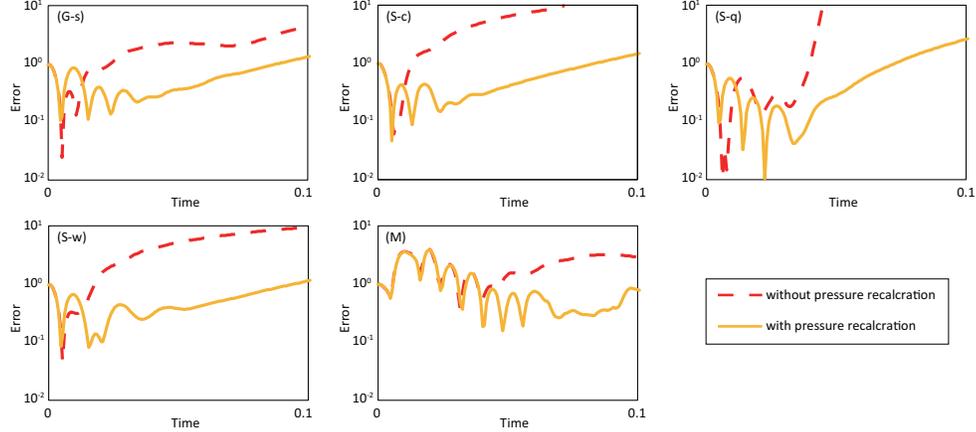}
	\caption{Time history of relative errors in space}
	\label{fig:TGvortex_error_time_history}
\end{figure}

Under the computational settings above, for the explicit particle method and that without the pressure recalculation \eqref{eq:disc_pres}, we compute the relative errors in space as:
\begin{equation}
\dfrac{\|\VelAppChar^{\timeindex{\tstep}}-\Vel^{\timeindex{\tstep}}\|_{\ell^2(\Dm)}}{\|\Vel^{\timeindex{\tstep}}\|_{\ell^2(\Dm)}},\qquad\dfrac{\|\PresAppave^{\timeindex{\tstep}}-\Pres^{\timeindex{\tstep}}\|_{\ell^2(\Dm)}}{\|\Pres^{\timeindex{\tstep}}\|_{\ell^2(\Dm)}}, 
\end{equation}
and the relative errors in space and time as:
\begin{equation}
\dfrac{\|\Vel-\VelAppChar\|_{\ell^2([0,T];\,\ell^2(\Dm))}}{\|\Vel\|_{\ell^2([0,T];\,\ell^2(\Dm))}},
\qquad\dfrac{\|\PresAppave-\Pres\|_{\ell^2([0,T];\,\ell^2(\Dm))}}{\|\Pres\|_{\ell^2([0,T];\,\ell^2(\Dm))}}.  
\end{equation}
Here, the norms are defined as
\begin{equation}
\|\phi^{\k}\|_{\ell^2(\Dm)} \deq \left(\sum_{\j=1}^{\N}\Pv{\j}|\phi^{\k}(\Pt{\j}{\k})|^2\right)^{1/2}, 
\end{equation}
\begin{equation}
\|\phi\|_{\ell^2([0,T];\,\ell^2(\Dm))}\deq \left(\sum_{\k=1}^{\tstepmax}\Dt\|\phi^{\k}\|_{\ell^2(\Dm)}\right)^{1/2}. 
\end{equation}
Moreover, $\PresAppave$ is defined by 
\begin{align}
\PresAppavek(\Pt{\i}{}) \deq \PresAppk(\Pt{\i}{})-\sum_{\j=1}^{\N}\Pv{\j}\PresAppk(\Pt{\j}{}),\quad \i=1,2,\dots,\N. 
\end{align}
Then, $\PresAppave$ satisfies the following condition: 
\begin{align}
\sum_{\i=1}^{\N}\Pv{\j}\PresAppave^k(\Pt{\j}{}) = 0,\quad k=1,2,\dots,K.
\end{align}
This condition corresponds to the integration condition of pressure: 
\begin{equation}
\int_{\Dm}\Pres(x,t)\dx = 0,\quad \forall t\in[0,T]. 
\end{equation}
Figure \ref{fig:TGvortex_error_time_history} shows time histories of the relative errors; in the figure, the vertical axes are plotted on the logarithmic scale. 
Table \ref{tab:TGvortex_error} lists the relative errors in space and time of velocity as well as pressure. 
In all the cases except for (M), the errors of pressure, which first oscillate, are considerably improved with the pressure recalculation compared with the case without the pressure recalculation. 
Moreover, the accuracy of velocity is enhanced by improving the accuracy of pressure.  
In the case of (M), the accuracy of pressure is not remarkably different between the cases with and without pressure recalculation; nevertheless, the accuracy of velocity is improved with pressure recalculation, which is clear from Table \ref{tab:TGvortex_error}.
Consequently, we use the method involving the pressure recalculation.

\begin{table}[h]
	\centering
	\caption{Relative errors in space and time for the cases with and without the pressure recalculation \eqref{eq:disc_pres}}
	\label{tab:TGvortex_error}
	\begin{tabular}{lrrr}
		& \multicolumn{3}{c}{(a) velocity} \\
		& without \eqref{eq:disc_pres} & with \eqref{eq:disc_pres} & without/with
		\\\hline
		(G-s) & $0.941$ & $0.022$ & $41.68$  \\
		(S-c) & $0.645$ & $0.030$ & $21.29$\\
		(S-q) & $164.784$ & $0.034$ & $4704.57$ \\
		(S-w) & $0.686$ & $0.028$ & $24.07$ \\
		(M) & $0.531$ & $0.034$ & $15.30$ 
	\end{tabular}
	\quad
	\begin{tabular}{rrr}
		\multicolumn{3}{c}{(b) pressure}\\
		without \eqref{eq:disc_pres} & with \eqref{eq:disc_pres} & without/with
		\\\hline
		$1.417$ & $0.520$  & $2.83$\\
		$4.892$ & $0.479$& $10.20$\\
		$2868.820$ & $0.572$& $5013.38$\\
		$3.515$ & $0.467$& $7.51$\\
		$2.209$ & $1.911$& $1.16$
	\end{tabular}
\end{table}

\begin{figure}[t]
	\centering
	\includegraphics[bb=0 0 108mm 75mm, width=60mm]{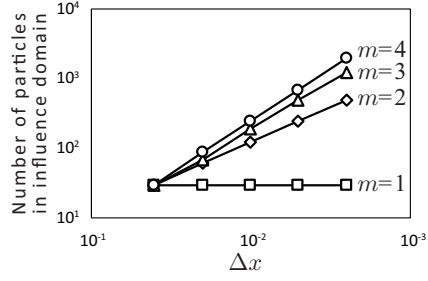}
	\caption{Graph of $\Dx$ versus average number of particles in the influence domain}
	\label{fig:ave_num_of_part_in_Bh}
\end{figure}

Next, we investigate the convergence of approximate solutions for the explicit particle method. 
We set the initial particle distributions using \eqref{TGvortex_ini_part_dist} with $\Dx=0.04, 0.02, 0.01, 0.005, 0.0025$. 
The particle volume set and reference weight functions are set as the same in the previous cases. 
For $m=1, 2, 3, 4$, the influence radius $\h$ is set as $\h=C_m \Dx^{1/m}$.  
Here, $C_m$ can be obtained by $C_m=3.1\times(0.04)^{1-1/m}$, which satisfies $h=3.1\times0.04$ when $\Dx=0.04$ for all $m$. 
We set $\Dt=\Dtmax$ and $\CompCof=2.5\Dx$. 
Then, by assuming that particles maintain distances proportional to $\Dx$, i.e., $r_{\min}=\order(\Dx)$, we have $\CompCof=\order(\Dx)$, $\h=\order(\Dx^{1/m})$, and $r_{\min}\h^{-2}=\order(\Dx^{1-2/m})$. 
Therefore, the conditions \eqref{cond:convergence} are satisfied when $\Dx\ra0$ in the case that $m>2$. 
It should be noted that although $m=1$ is often used in practical computing because the average number of particles in the influence domain increases exponentially when $m>1$, as shown in Figure \ref{fig:ave_num_of_part_in_Bh}, we must set $m>1$ to conduct simulations with convergence.

\begin{figure}[!t]
	\centering
	\includegraphics[bb=0 0 181mm 139mm, width=150mm]{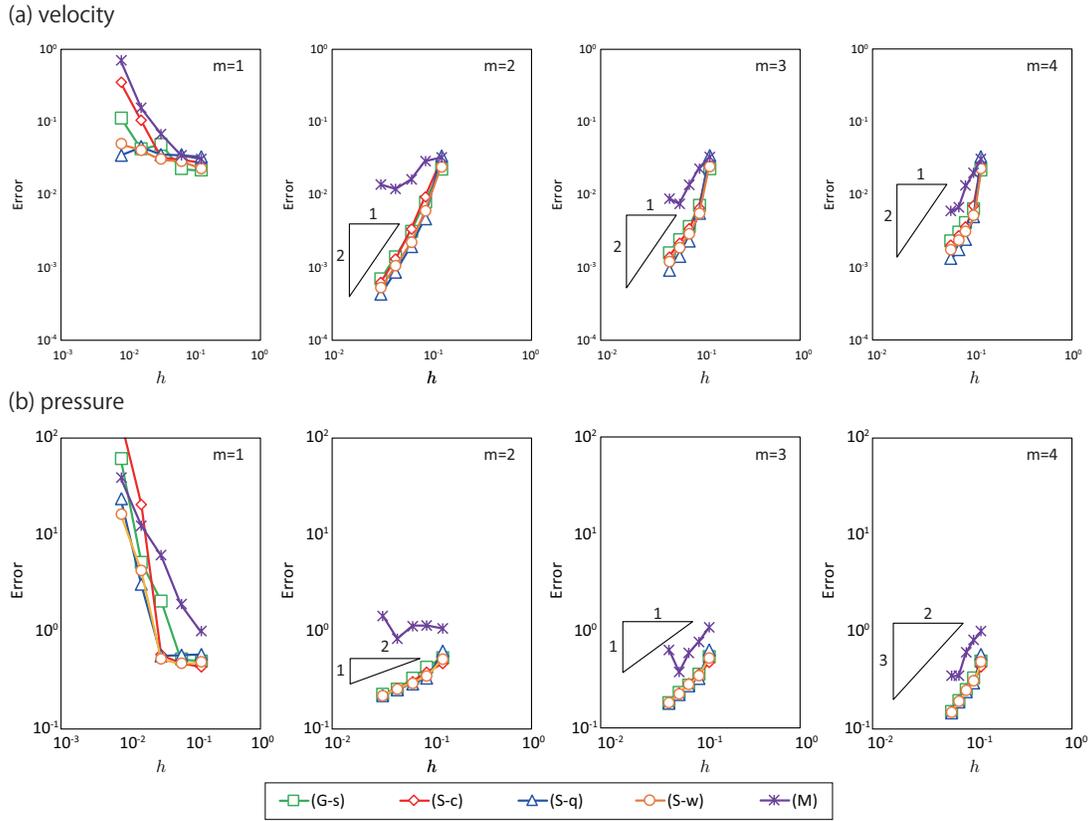}
	\caption{Graphs of relative errors versus the influence radius}
	\label{fig:TGvortex_Re0010_error_convergence}
\end{figure}

Under these conditions, we compute the relative errors in space and time. 
Figure \ref{fig:TGvortex_Re0010_error_convergence} shows the double logarithmic graph of the relative errors versus the influence radius $h$. 
Here, the slopes of the hypotenuse of the triangle in Figure \ref{fig:TGvortex_Re0010_error_convergence} ($m\geq2$) are $\mathcal{O}(h^{2})$ for (a) velocity and $\mathcal{O}(h^{(m-1)/2})$ for (b) pressure. 
Table \ref{tab:TGvortex_Re0010_error_convergence} lists the convergence rates of velocity and pressure obtained for $\Dx=0.005$ and $0.0025$. 
From Figure \ref{fig:TGvortex_Re0010_error_convergence} and Table \ref{tab:TGvortex_Re0010_error_convergence}, we can confirm that the convergence orders of velocity and pressure with respect to the influence radius $\h$ are of the second order and $(m-1)/2$th order for $m\geq2$, respectively, except in case (M). 
This is because the approximate solution did not converge in the case of the approximate operators of MPS; this might be attributed to the fact that sufficient conditions of convergence were derived under assumptions of sufficiently smooth and bounded weight functions for the truncation error estimates. 

\begin{table}[!h]
	\centering
	\caption{Convergence rates of errors obtained when $\Dx=0.005$ and $0.0025$}
	\label{tab:TGvortex_Re0010_error_convergence}
	\centering
	\begin{tabular}{lrrrr}
		& \multicolumn{4}{c}{(a) velocity}\\
		& $m=1$ & $m=2$& $m=3$ & $m=4$ 
		\\\hline
		(G-s) & -1.43  & 2.13  & 1.84 & 1.71 \\
		(S-c) & -1.73  & 2.18  & 1.84 & 1.69 \\
		(S-q) & 0.37  & 2.14  & 1.84  & 1.64 \\
		(S-w) & -0.29   & 2.07  & 1.82  & 1.69 \\
		(M) & -2.18  & -0.42 & -0.78 & 0.79 
	\end{tabular}~
	\begin{tabular}{rrrr}
		\multicolumn{4}{c}{(b) pressure}\\
		$m=1$ & $m=2$& $m=3$ & $m=4$
		\\\hline
		-3.54  & 0.35  & 1.07 & 1.56\\
		-2.89  & 0.50  & 0.98 & 1.55\\
		-2.94  & 0.46   & 0.93  &  1.49  \\
		-1.92  & 0.48  & 0.95  & 1.52   \\
		-1.65 & -1.65 & -2.29 & 0.08
	\end{tabular}
\end{table}

\section{Approaches for reducing truncation errors of the generalized approximate operators}
\label{sec:improve_operator}
In order to conduct more accurate simulations, we improve the accuracy of the generalized approximate operators by considering an optimization problem for weight functions derived based on their truncation error estimates.
Moreover, the accurate results of the explicit particle methods (i.e., with an optimal weight function included) are confirmed by numerical truncation errors and numerical errors of cavity flow. 

\subsection{Optimization problem for weight functions}
\label{subsec:optimal_weight_func}
We derive an optimization problem of truncation errors for weight functions and its solutions. 
As discussed in Section \ref{subsec:truncation_error}, the truncation error of the generalized approximate Laplace operator is estimated by 
\begin{equation}
|\lap\f_\i-\LapApp{}\f_\i| \leq |\widetilde{E}_{\i}| + |\widehat{E}_{\i}|,
\label{terror:lap}
\end{equation}
where
\begin{equation}
	\widetilde{E}_{\i} = \dfrac{4\dim}{C_{2}(\wl)\h^2}\int_{\DmH}R_{4,\i}(y;\f)\wlh(|y-\Pt{\i}{}|)\dy,
	\end{equation}
	\begin{equation}
	\widehat{E}_{\i} = \dfrac{2\dim}{C_{2}(\wl)\h^2}\brm{\int_{\DmH}\{\f(y)-\f_\i\}\wlh(|y-\Pt{\i}{}|)\dy
		-\sum_{\j=1}^\N\Pv{\j} (\f_\j-\f_\i)\wlh(|\Pt{\j}{}-\Pt{\i}{}|)}. 
	\end{equation}
We can estimate $|\widetilde{E}_{\i}|$ using $\order(\h^2)$ independent of particle distributions. 
Thus, we can estimate that $|\widehat{E}_{\i}|$ represents an error based on disturbances of the particle distribution. 
In practical computing, it is rare for the particle distribution to become sufficiently uniform in each time step; hence, we aim to reduce the error $|\widehat{E}_{\i}|$. 
In Section \ref{subsec:truncation_error}, we estimated $|\widehat{E}_{\i}|$ as
\begin{equation}
|\widehat{E}_{\i}| =\dfrac{2\dim}{C_{2}(\wl)}\brm{\int_{\dRd}\brs{\wl(|y|)+2\BrA{\deriv{1}{r}{}\wl(|y|)}}\dy}\frac{\delta_{\infty}}{\h^2}|\f|_{C^1(\DmHo)}+\order(\delta_{\infty}^2\h^{-3})
\end{equation}
under the condition $\w\in\FsWeightFunc\cap C^1([0,\infty))$. 
Therefore, by using this term with respect to the weight function in $|\widehat{E}_{\i}|$ as the objective function, we define the following optimization problem for the weight functions: 
\begin{align}
\left\|
\begin{array}{ll}
\mbox{minimize}&\displaystyle F(\w) =\brm{\int_{\dRd}\brs{\w(|y|)+2\BrA{\deriv{1}{r}{}\w(|y|)}}\dy}\brm{\int_{\dRd}|y|^2|\w(|y|)|\dy}^{-1}
\smallskip\\
\mbox{subject to}& \mbox{$\w$ satisfies $\w\in\FsWeightFunc\cap C^1([0,\infty))$}.
\end{array}
\right.
\label{weight:optimization_problem}
\end{align}
In order to reduce the computational complexity, we consider solving the optimization problem within a range where the reference weight function transforms into a polynomial function. 
We give the reference weight function as the $n$th polynomial function: 
	\begin{equation}
	w(r)=\sum_{k=0}^{n}a_kr^k, \qquad a_0,a_1,\dots,a_n\in\dR.
	\label{weight_poly}
	\end{equation}
	Because the condition $\w\in\FsWeightFunc\cap C^1([0,\infty))$ yields
	\begin{equation}
	w(0)>0,\qquad w(1)=0, \qquad \dfrac{d}{dr}w(1)=0,
	\end{equation}
	we have the conditions of the coefficients in \eqref{weight_poly}:
	\begin{equation}
	a_0>0,\qquad \sum_{k=0}^n a_k=0,\qquad \sum_{k=1}^n ka_k=0. 
	\label{cond:poly_coef}
	\end{equation}
	Therefore, in the case of a quadratic polynomial $n=2$, the solution of \eqref{weight:optimization_problem} is the spike function \eqref{weight_spike}. 
	When $n\geq3$, we consider that the additional condition $\mbox{minimizes}~F(\w)$ because we calculate $F(\w)$ as 
	\begin{align}
	F(\w)
	&=\brm{\int_{0}^1r^{\dim-1}\brs{\w(r)+2\BrA{\deriv{1}{r}{}\w(r)}}dr}\brm{\int_{0}^1r^{\dim+1}|\w(r)|dr}^{-1}\nn\\
	&=\brm{\int_{0}^1\brs{\sum_{k=0}^{n}a_kr^{k+\dim-1}+2\BrA{\sum_{k=1}^{n}ka_kr^{k+\dim-2}}}dr}\brm{\int_{0}^1\sum_{k=0}^{n}a_kr^{k+\dim+1}dr}^{-1}\nn\\
	&=\brs{\sum_{k=0}^{n}\dfrac{a_k}{k+\dim}+2\int_{0}^1\BrA{\sum_{k=1}^{n}ka_kr^{k+\dim-2}}dr}\brs{\sum_{k=0}^{n}\dfrac{a_k}{k+\dim+2}}^{-1}
	\end{align}
	By solving the minimization problem of $F(w)$ for the coefficients of polynomial functions under the constraint conditions \eqref{cond:poly_coef}, we can obtain optimal weight functions for $n\geq3$. 
	However, because the optimal weight functions with $n\geq3$ depend on the spatial dimension, we use the quadratic spike function \eqref{weight_spike} to avoid such spatial dimension dependency in the subsequent numerical experiments. 

\subsection{Numerical results of truncation errors}
\label{subsec:numerical_truncation_err}
In order to verify the analytical discussions presented in the previous section, we compute the numerical truncation errors of approximate Laplace operators when the disturbances of the particle distribution are changed. 
Furthermore, the test function is set as $v(\coord{x}{1},\coord{x}{2})=\sin(2\pi(\coord{x}{1}+\coord{x}{2}))$. 
The domain is set as $\Dm=(0,1)\times(0,1)$. 
Let $H=3\times2^{-4}$. 
Then, the particle distribution is set as
\begin{align}
\PtSet=\{((i-1/2+\Pertab_{\i\j}^{(1)}/2)\Dx, (j-1/2+\Pertab_{\i\j}^{(2)}/2) \Dx)\in\DmH ;~i,j\in\dZ\}.
\label{particle_distribution_random_01}
\end{align}
Here, $\Dx=2^{-4}$ and $\Pertab_{\i\j}^{(k)}\,(k=1,2)$ is a random number satisfying $|\Pertab_{\i\j}|\leq\PertabMax\,(0\leq \PertabMax< 1)$. 
Figure \ref{figpd} shows examples of the particle distributions with perturbation $\PertabMax=0, 0.25, 0.5$ in $\Dm$. 
\newcommand{\figuresizeparticledistribution}{130mm}
\begin{figure}[t]
	\centering
	\includegraphics[bb=0 0 435.0mm 141.3mm, width=\figuresizeparticledistribution]{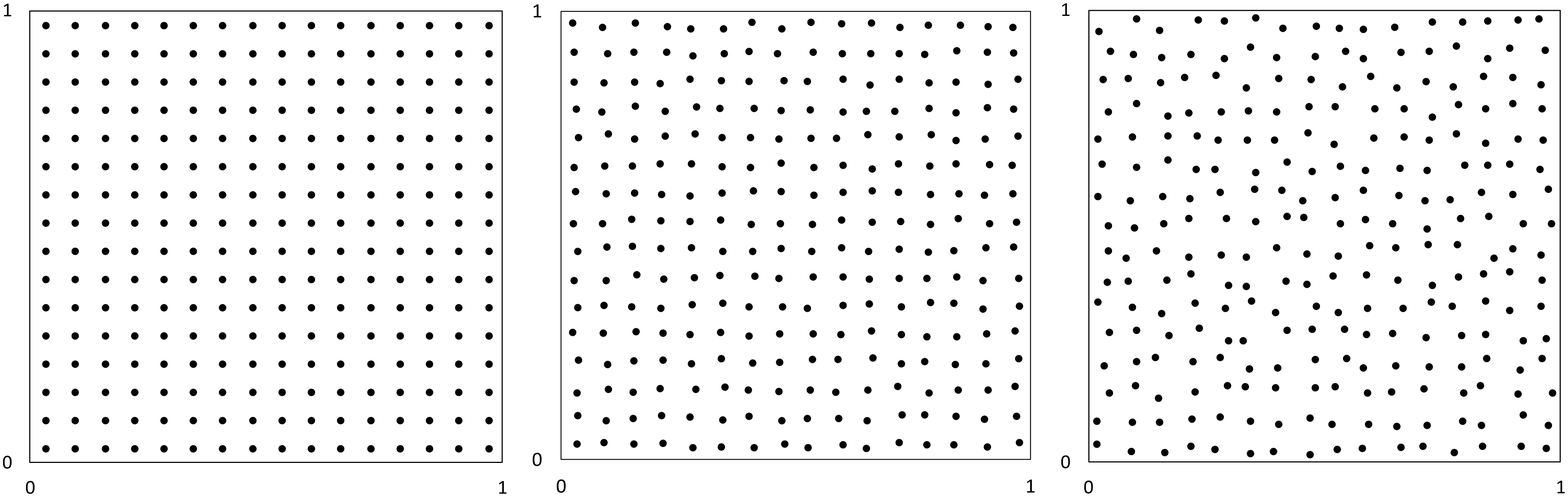}
	\caption{Particle distributions with $\PertabMax=0$ (left), $\PertabMax=0.25$ (center), and $\PertabMax=0.5$ (right). }
	\label{figpd}
\end{figure}
It should be noted that this particle distribution becomes a square lattice if $\PertabMax=0$. 
The particle volume set $\PvSet=\{\Pv{\i}\}_{\i=1}^\N$ is determined by \eqref{pt_volume_eg}. 
The influence radius is set as $\h=2.1\Dx, 2.6\Dx, 3.1\Dx$. 
We consider the following four reference weight functions: 
\begin{enumerate}
	\setlength{\itemindent}{3.5mm}
	\setlength{\parskip}{2ex}
	\setlength{\itemsep}{0ex}
	\item[(G-s)] $\wl$ is set as the quadratic spike function \eqref{weight_spike};
	\item[(S-c)] $\wl$ is set as 
	\begin{equation}
	\wl(r)=-\dfrac{1}{r}\WSPHDriv(r), 
	\label{weight_SPH_lap}
	\end{equation}
	\hspace{2ex}where $\ws$ is the cubic B-spline \eqref{w_SPH_cubic}. 
	\item[(S-q)] $\wl$ is set as \eqref{weight_SPH_lap} where $\ws$ is the quintic B-spline \eqref{w_SPH_quintic};
	\item[(S-w)] $\wl$ is set as \eqref{weight_SPH_lap} where $\ws$ is the quintic Wendland function \eqref{w_SPH_wendland}.
\end{enumerate}

\begin{figure}[t!]
	\centering
	\includegraphics[bb=0 0 438.0mm 209mm, width=156mm]{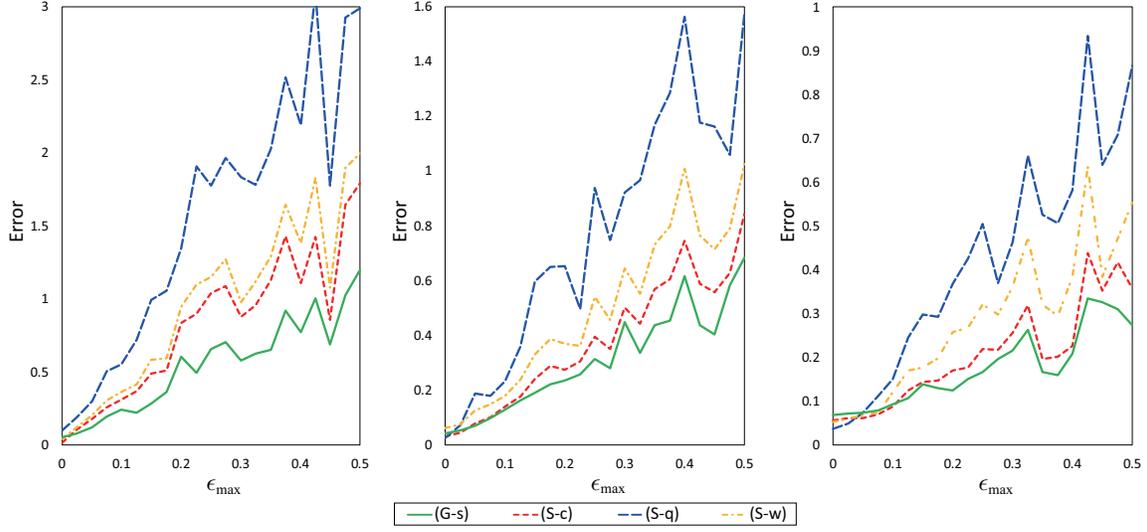}
	\caption{Graphs of the relative truncation errors of approximate Laplace operators versus perturbation $\PertabMax$ with $\h=2.1\Dx$ (left), $\h=2.6\Dx$ (center), and $\h=3.1\Dx$ (right). }
	\label{fig_terror_lap_perturbation}
\end{figure}

Figure \ref{fig_terror_lap_perturbation} shows the graphs for the relative truncation error
\begin{align}
\frac{\displaystyle \max_{\Pt{\i}{}\in\Dm}|\Delta \f_\i-\LapApp{} \f_\i|}{\displaystyle \max_{\Pt{\i}{}\in\Dm}|\Delta \f_\i|}
\end{align}
versus the perturbation $\PertabMax$ when $\h=2.1\Dx, 2.6\Dx, 3.1\Dx$. 
Table \ref{table_terror_lap_perturbation} lists the relative truncation errors with $\h=2.1\Dx, 2.6\Dx, 3.1\Dx$ and $\PertabMax=0.0, 0.25, 0.5$. 
From Figure \ref{fig_terror_lap_perturbation} and Table \ref{table_terror_lap_perturbation}, we can confirm that the truncation errors increase as perturbation $\PertabMax$ increases and influence radius rate $\h$ decreases. 
In all the cases, though the truncation error of the generalized approximate Laplace operator with the spike function is larger than that of the conventional Laplace operators for uniform particle distributions ($\PertabMax=0$), the truncation error becomes smaller for general particle distributions ($\PertabMax>0$). 
Therefore, we confirmed that truncation errors can be effectively reduced for general particle distributions using the generalized Laplace operator with the spike function. 
Later, in Section \ref{subsec:Driven_cavity_flow}, we confirm whether the generalized approximate operators with the spike function are also valid for fluid simulations. 

\begin{table}[h!t]
	\caption{Truncation errors of approximate Laplace operators with $\h=2.1\Dx, 2.6\Dx, 3.1\Dx$ and $\Pertab=0.0, 0.25, 0.5$. }
	\label{table_terror_lap_perturbation}
	\centering
	\newcommand{\tablespaceA}{~~}
	\begin{tabular}{l@{\quad}r@{\tablespaceA}r@{\tablespaceA}r}
		& \multicolumn{3}{c}{$\h=2.1\Dx$} \\
		& $\PertabMax=0$ & $\PertabMax=0.25$ &  $\PertabMax=0.5$ 
		\\\hline
		(G-s) & 0.0532 & 0.6609 & 1.2014\\
		(S-c)  & 0.0191 & 1.0407 & 1.7955\\
		(S-q) & 0.0994 & 1.7798 & 2.9894\\
		(S-w)  & 0.0447 & 1.1538 & 1.9976
	\end{tabular}
	\quad
	\begin{tabular}{l@{\quad}r@{\tablespaceA}r@{\tablespaceA}r}
		& \multicolumn{3}{c}{$\h=2.6\Dx$} \\
		& $\PertabMax=0$ & $\PertabMax=0.25$ &  $\PertabMax=0.5$ 
		\\\hline
		(G-s) & 0.0409 & 0.3143 & 0.6837\\
		(S-c)  & 0.0306 & 0.3934 & 0.8428\\
		(S-q) & 0.0296 & 0.9434 & 1.5808\\
		(S-w)  & 0.0607 & 0.5408 & 1.0244
	\end{tabular}
	\quad
	\begin{tabular}{l@{\quad}r@{\tablespaceA}r@{\tablespaceA}r}
		& \multicolumn{3}{c}{$\h=3.1\Dx$} \\
		& $\PertabMax=0$ & $\PertabMax=0.25$ &  $\PertabMax=0.5$ 
		\\\hline
		(G-s) & 0.0695 & 0.1673 & 0.2731\\
		(S-c)  & 0.0567 & 0.2205 & 0.3589\\
		(S-q) & 0.0383 & 0.5048 & 0.8656\\
		(S-w)  & 0.0529 & 0.3210 & 0.5537
	\end{tabular}
\end{table}

\subsection{Numerical results of driven cavity flow}
\label{subsec:Driven_cavity_flow}
In order to investigate whether the generalized approximate operators with the optimal weight function are also effective for a flow problem, we apply the explicit particle method to a driven cavity and compare errors in the cases when five pairs of weight functions are used. 
The driven cavity flow is a viscous flow problem in a rectangular domain with Dirichlet boundary conditions. 
One side of the boundaries flows in a tangential direction, while the other sides are wall boundaries. 
In the case of a square domain, the driven cavity flow can be denoted solely on the basis of the Reynolds number $Re=LU/\nu$, where $L$ and $U$ are the length of one side of the domain and velocity on the driven boundary, respectively. 
Hereafter, we consider $Re=100$ and $1000$. 

We consider the driven cavity flow in the square domain $\Dm=(0,1)\times(0,1)$. 
We denote the velocity as $\Vel=(\Vel_1, \Vel_2)^{T}$. 
The initial conditions are given by
\begin{align}
(\Vel_1, \Vel_2)^{T} = (0,0)^{T},\qquad (x,y)\in \Dm,~t=0 
\end{align}
while the boundary conditions are given by
\begin{align}
(\Vel_1, \Vel_2)^{T} = 
\begin{cases}
(1,0)^{T}\qquad& \mbox{in~}\brm{(x,t)\in\bd\times(0,\tmax);\,\coord{x}{2}= 1},
\\
(0,0)^{T}\qquad &\mbox{in~}\brm{(x,t)\in\bd\times(0,\tmax);\,\coord{x}{2}< 1}. 
\end{cases}
\end{align}
Furthermore, zero gravity is assumed ($f=0$). 

We set the parameters as follows. 
Set $H=0.1$, $\Veliniex=0$ in $\DmH$ and $\Velbdex=(1,0)  \mbox{~in~}\{(x,t)\in\DmH\times(0,\tmax);\,\coord{x}{2}\geq1\}$, $=(0,0)  \mbox{~in~}\{(x,t)\in\DmH\times(0,\tmax);\,\coord{x}{2}<1\}$ . 
The initial particle distribution $\PtSet^{\timeindex{0}}$ is set as a square lattice with spacing $\Dx$. 
Here, $\Dx$ is set as $\Dx=0.005$ and $\Dx=0.0025$ when $Re=100$ and $Re=1000$, respectively. 
Note that the particles are distributed in $\DmH=(-H,1+H)\times(-H,1+H)$, and the particle distributions outside of the wall boundary correspond to well-known dummy particles \cite{shao2003incompressible}. 
We consider the same five pairs of approximate operators used in Section \ref{subsec:taylor_green}. 
We consider the same three cases of influence radii as in the numerical experiments in Section \ref{subsec:numerical_truncation_err}: (a) $\h=2.1\Dx$; (b) $\h=2.6\Dx$; (c) $\h=3.1\Dx$. 
We set $\CompCof=0.1$ and $\Dt=\Dtmax$. 
Under these conditions, we compute the two-dimensional driven cavity flow and compare velocity profiles in the vertical direction on the lines $x=0.5$ with the reference solutions, which are the numerical results of the higher-order finite difference method by Ghia et al.\ \cite{ghia1982high}. 

\begin{figure}[th]
	\newcommand{\figsizeA}{0.625}
	\begin{center}
		\includegraphics[bb=0 0 344.7mm 207.1mm, width=160mm]{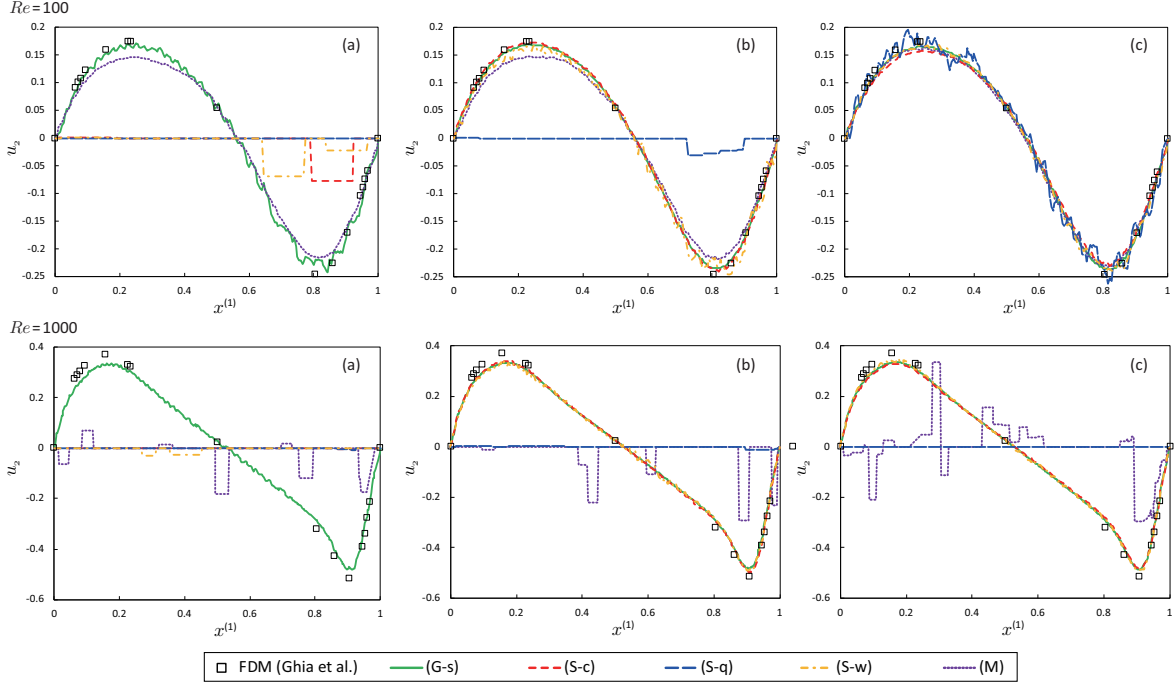}
		\caption{Velocity profiles of two-dimensional driven cavity flow at $Re = 100\,(\Dx=0.005)$ and $1000\,(\Dx=0.0025)$. (a) $\h=2.1\Dx$; (b) $\h=2.6\Dx$; (c) $\h=3.1\Dx$. 
		}
		\label{fig:cavity}
	\end{center}
\end{figure}

Figure \ref{fig:cavity} shows the velocity profiles of the two-dimensional driven cavity flow at $Re = 100$ and $Re=1000$.
The boxes in Figure \ref{fig:cavity} show the vertical velocity $\Vel_2^{\rm\,FDM}(x^{\rm\,FDM}_j)$ of the observation point $x^{\rm\,FDM}_j$ in the results of Ghia et al. \cite{ghia1982high}. 
Table \ref{tab:cavity_L2_error_Re0100} lists the errors of velocity measured using the following discrete $L^2$ norm in space: 
\begin{align}
\brs{\dfrac{\ds\sum_{j=1}^{M} \Delta \Pt{\j}{\rm\,FDM} \left|\Vel_2(x_{i_j})-\Vel_2^{\rm\,FDM}(x^{\rm\,FDM}_j)\right|^2}{\ds\sum_{j=1}^{M}\Delta \Pt{\j}{\rm\,FDM}\left|\Vel_2^{\rm\,FDM}(x^{\rm\,FDM}_j)\right|^2}}^{1/2},
\end{align}
where $\Delta \Pt{\j}{\rm\,FDM}=\Pt{\j}{\rm\,FDM}-x_{j-1}^{\rm\,FDM}$, and   
\begin{align}
i_j = \argmin_{k}|\Pt{\j}{\rm\,FDM}-x_k|. 
\end{align}
From the Figure \ref{fig:cavity} and Table \ref{tab:cavity_L2_error_Re0100}, it is clear that the velocity becomes stable as the influence radius increases. 
In particular, because case (G-S) has a solution even in the case of (a) $\h=2.1\Dx$, the explicit particle method using the generalized approximate operators with the spike weight function \eqref{weight_spike} is more robust to the influence radius than that with other weight functions. 
This result is consistent with that of the generalized approximate operators with the spike weight function \eqref{weight_spike} being more accurate for non-uniform particle distributions than other operators for the truncation error estimates, as discussed in Sections \ref{subsec:optimal_weight_func}--\ref{subsec:numerical_truncation_err}. 
Therefore, we confirm that the generalized approximate operators with the spike function \eqref{weight_spike} are also effective for a flow problem. 

\newcommand{\tablespaceCavity}{\hspace{3ex}}
\begin{table}[!htbp]
	\begin{center}
		\caption{Relative errors between the reference solutions (FDM) and particle methods.}
		\label{tab:cavity_L2_error_Re0100}
		\begin{tabular}{lc@{\tablespaceCavity}c@{\tablespaceCavity}c@{\hspace{6ex}}c@{\tablespaceCavity}c@{\tablespaceCavity}c}
			&\multicolumn{3}{c}{$Re=100$}&\multicolumn{3}{c}{$Re=1000$}\\
			& (a) & (b) & (c) & (a) & (b) & (c)
			\\\hline
			(G-s) & 0.0667 & 0.0484 & 0.0486 & 0.0975 & 0.0904 & 0.0877\\
			(S-c) & 0.7546 & 0.0513 & 0.0914 & 0.9987 & 0.0786 & 0.1068 \\
			(S-q) & 1.0000 & 0.9206 & 0.0816 & 0.9987 & 0.9916 & 0.9998\\
			(S-w) & 0.9809  & 0.1206  & 0.0592 & 1.0002 & 0.0914  & 0.0835\\
			(M) & 0.1326 & 0.1299 & 0.0698 & 1.0757 & 0.9397 & 1.2913
		\end{tabular}
	\end{center}
\end{table}

\section{Application for incompressible viscous flow problems under free surface effects}
\label{sec:application_free_surface}
In order to confirm that the explicit particle method is applicable to realistic problems, we develop the explicit particle method for flow problems under free surface effects.
We introduce modifications for pressure evaluation and pressure gradient to avoid clustering of particles in and around a free surface.  
Moreover, we apply the modified explicit particle method to a dam break flow and compare the numerical results with experimental results.

\subsection{Treatment of free surface}
\label{subsec:treat_free_surface}
Because particles around the free surface do not have a sufficient number of particles in their influence domain, approximate operators on these particles do not behave appropriately. 
In particular, particles around the free surface come close or collide with each other in the case of the original scheme. 
For this reason, the tentative densities on particles around the free surface are evaluated to be considerably lower than that an inner particle.  
Consequently, the tentative pressure on these particles become negative per \eqref{eq:disc_pres_temp}; then, retraction forces are experienced owing to the pressure gradient in \eqref{eq:disc_vel}. 
In order to solve this problem, we have to modify evaluations of the tentative density and tentative pressure. 
Thus, we modify \eqref{eq:disc_pres_temp} and \eqref{eq:disc_vel}. 
In order to avoid obtaining negative pressures, we modify \eqref{eq:disc_pres_temp} as
\begin{equation}
\PresAppP{\i}{\k+1} =  \max\brm{\dfrac{\rho}{\CompCof^{2}}\brs{\dfrac{1}{C_{0,h}(\w)}\sum_{\j=1}^{\N}\Pv{\j}\wh(|\PtP{\j}{\k+1}-\PtP{\i}{\k+1}|)-1}, 0}.
\label{eq:disc_pres_temp_modify}
\end{equation}
Moreover, when the original pressure gradient is used in \eqref{eq:disc_vel}, a non-physical force develops in the tangential direction of the free surface because of the lack of a sufficient number of particles in and around the free surface.
Therefore, we modify \eqref{eq:disc_vel} as
\begin{equation}
\left\{
\begin{array}{@{\,}r@{\,}c@{\,}ll}
\dfrac{\VelApp{\i}{\k+1}-\VelAppP{\i}{\k+1}}{\Dt} &=& -\dfrac{1}{\rho}\GradAppPlus{\k+1}\PresApp{\i}{\k+1},\qquad&\i\in\IndexSet{\Dm}{\k},
\\
\VelApp{\i}{\k+1}&=&\Velbdex(\Pt{\i}{\k+1},\tkp),\qquad&\i\in\IndexSet{\bdh}{\k}.
\end{array}
\right.
\label{eq:disc_vel_modify}
\end{equation}
Only with the modifications above, we attain a stable and accurate simulation of a dam break flow in the next section; however, we observe strange particle motions around the free surface. 
Therefore, we add the collision methods \cite{shakibaeinia2012mps} used in E-MPS, which modify the particle distributions to maintain their momentum, and this was confirmed to solve the problem. 

\subsection{Dam break flow}
\label{subsec:dam_break}
The dam break flow is a flow problem in which a water column on one side of a tank collapses because of gravity.   
Because a considerable amount of experimental data, including flow tip speeds, wave height history, and wall pressure distributions, have been collected in the literature \cite{lobovsky2014experimental,hu2004cip}, changes in free surface geometry and pressure distributions in the numerical results can be confirmed.   

\begin{figure}[t]
	\begin{center}
		\includegraphics[bb=0 0 852.7mm 277.6mm, scale=0.19]{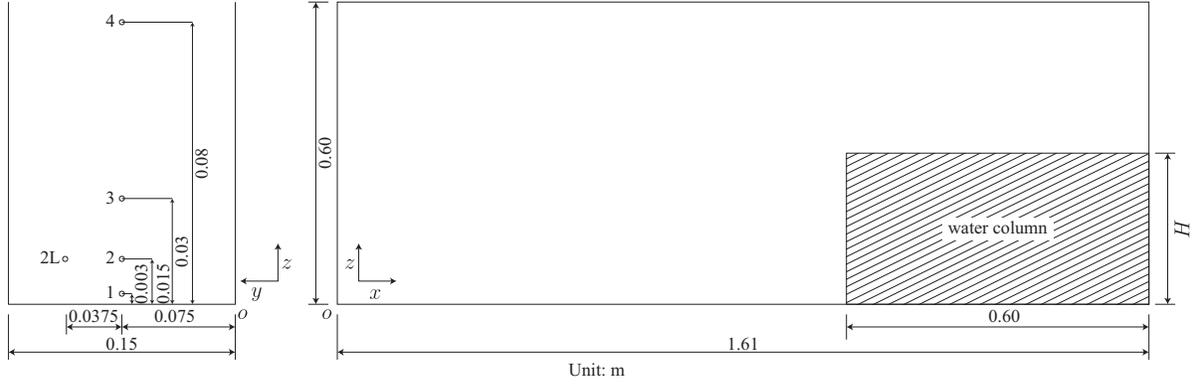}
		\caption{Computational model of three-dimensional dam break flow.}
		\label{dambreakmodel01}
	\end{center}
\end{figure}

We consider the hydraulic experiment by Lobovsk{\`y} et al.\ \cite{lobovsky2014experimental} as shown in Figure \ref{dambreakmodel01}. 
In this experiment, five pressure sensors are set on the opposite side of the water column. 
As shown in the left part of Figure \ref{dambreakmodel01}, the five pressure sensors labeled as 1, 2, 2L, 3, and 4.
In particular, their coordinates are $(0,0.075,0.003)$, $(0,0.075,0.015)$,  $(0,0.0375,0.015)$, $(0,0.075,0.03)$, and $(0,0.075,0.08)$, respectively, from the origin $o$. 
The height of the water column $H_{\rm dam}$ is $0.3$ m or $0.6$ m. 

We set the end time $T$ as $T=1.3$. 
Furthermore, we set the remaining parameters as follows. 
The initial particle distribution in the flow domain is set as a cubic lattice with spacing $\Dx=3.0\times10^{-3}$ m in the water column. 
Moreover, we set particles on a cubic lattice on the outer domain whose distances from the wall are less than $H=5.2\Dx$. 
Note that the particle distributions outside of the wall boundary correspond to well-known dummy particles \cite{shao2003incompressible}. 
The velocity of the particles outside the domain are set as zero. 
Then, we set $\CompCof=0.05$, $\h=2.6\Dx$, and $\Dt=\Dtmax$. 
Under these conditions, we compute the dam break flow and compare the pressure at the sensors.
Here, the pressures at the sensors are computed using the pressure of the nearest particle on the wall boundary from these sensors, i.e., the numerical pressure of Sensor $l$ at $t=t^k$ is computed as
\begin{align}
P_{l}(t^k)  = p(x_{i_l}^k), \qquad i_l  = \argmin_{k} |X_l - x_{k}|, 
\end{align}
where $X_l$ is the position of Sensor $l$. 

\begin{figure}[t]
	\begin{center}
		\includegraphics[bb=0 0 418.5mm 231.0mm, scale=0.38]{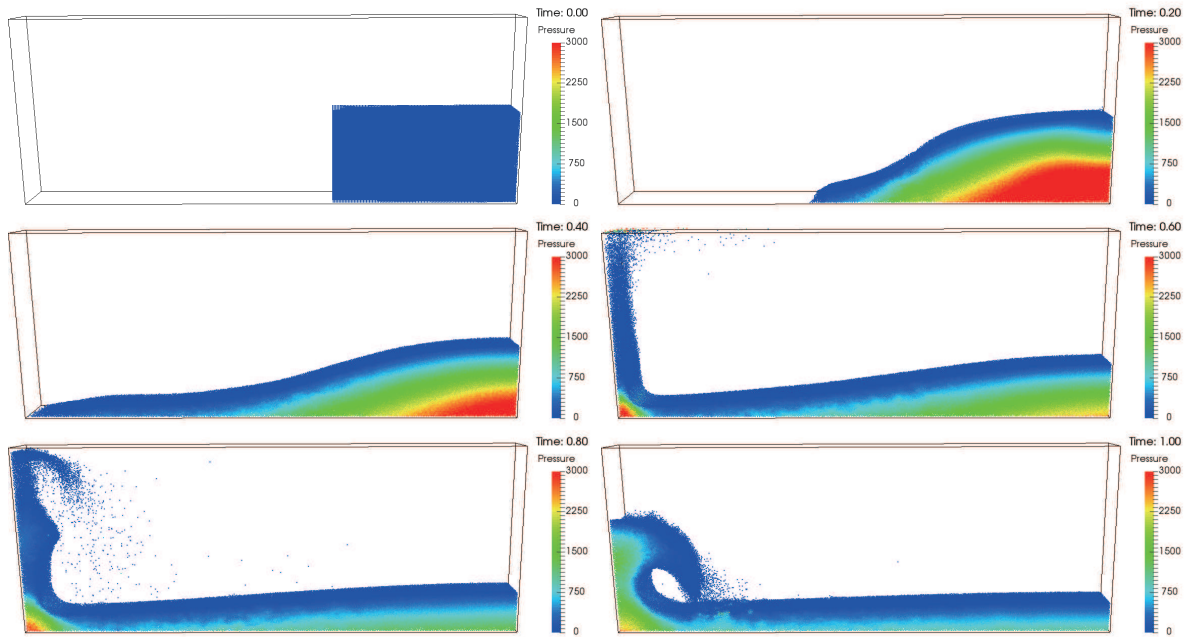}
		\caption{Pressure distributions for three-dimensional dam break flow ($H_{\rm dam}=0.3$ m).}
		\label{fig:dambreak_pressure_contor}
	\end{center}
\end{figure}

\begin{figure}[th]
	\includegraphics[bb=0 0 650.5mm 681.7mm, scale=0.25]{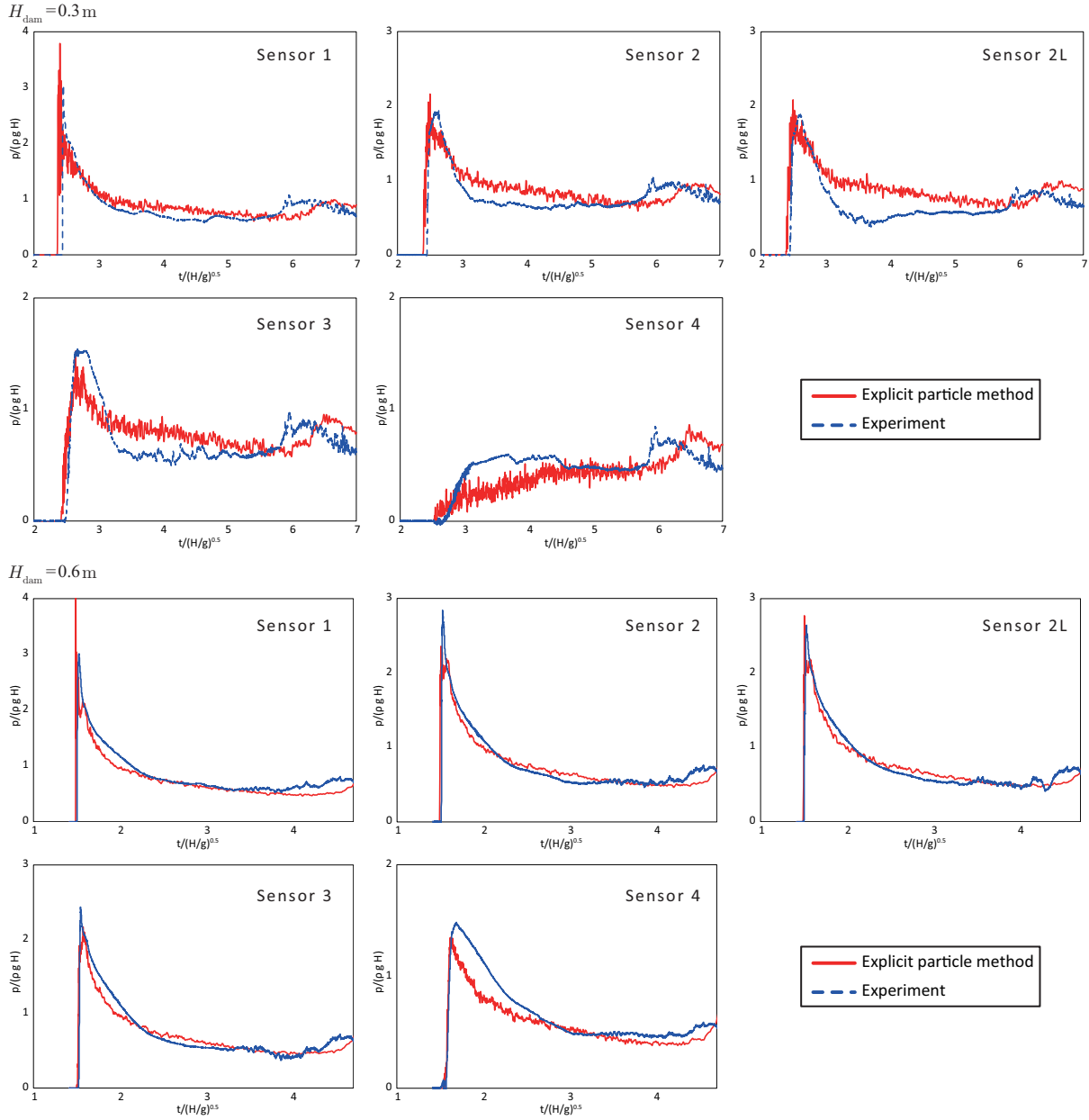}
	\caption{Pressure history at sensors for dam break flow.}
	\label{fig:dambreak_pres_histry}
\end{figure}

Figure \ref{fig:dambreak_pressure_contor} shows the pressure distributions of the explicit particle method when $H_{\rm dam}=0.3$. 
Figure \ref{fig:dambreak_pres_histry} shows pressure histories of the experimental and numerical results at each sensor when $H_{\rm dam}=0.3$\,m and $0.6$\,m. 
Table \ref{tab:dambreak_errors} lists relative errors of pressure in a discrete $L^2$ norm in time as 
\begin{align}
\dfrac{\sqrt{\sum_{k=1}^{K}\Delta t^{k}\left|P_{l}(t^k)-P^{\rm\,ex}_{l}(t^k)\right|^2}}{\sqrt{\sum_{k=1}^{K}\Delta t^{k}\left|P^{\rm\,ex}_{l}(t^k)\right|^2}}
\end{align}
for Sensor $l$. 
Here, $P^{\rm\,ex}_{l}(t^k)$ is the observed pressure for Sensor $l$ at $t=t^k$. 
From Figure \ref{fig:dambreak_pressure_contor}, we can observe smooth pressure distributions. 
Moreover, from Figure \ref{fig:dambreak_pres_histry} and Table \ref{tab:dambreak_errors}, we can obtain the numerical results based on the experiment results. 
These numerical results show that the explicit particle method is applicable for flow problems under free surface effects. 

\begin{table}[h]
	\begin{center}
		\caption{Relative errors between experiments and particle methods at sensors.}
		\label{tab:dambreak_errors}
		\begin{tabular}{rrrrrr}
			$H_{\rm dam}$& 1 & 2 & 2L & 3 & 4
			\\\hline
			0.3 & 0.3829 & 0.2374 & 0.2363 & 0.2161 & 0.1834 \\
			0.6 & 0.2939 & 0.2307 & 0.2275 & 0.1953 & 0.1998 
		\end{tabular}
	\end{center}
\end{table}

\section{Conclusion}
\label{sec:conclusion}
We conducted a convergence study for an explicit particle method for the incompressible Navier--Stokes equations. 
The explicit particle method is based on a penalty problem of the incompressible Navier--Stokes equations, which was derived using the mathematical discretization procedure. 
Moreover, the explicit particle method uses generalized approximate operators, which was introduced as a wider class of approximate operators than those used in SPH and MPS for spatial discretization. 
By referring to the convergence orders of the penalty problem and orders of the residual appearing in the derivation process as well as truncation errors of the generalized approximate operators, we conjectured sufficient conditions of convergence for the explicit particle method. 
The convergence with these sufficient conditions was confirmed using numerical results of the Taylor--Green vortex; in particular, these numerical convergence orders of velocity and pressure with respect to the influence radius $h$ were $\order(\h^2)$ and $\order(\h^{(m-1)/2})$ with $m\geq2$, respectively, where $m$ is a parameter determining the ratio of increase of neighbor particles in influence.

Next, we optimized the reference weight functions considering the decreasing truncation errors of the generalized approximate operators for non-uniform particle distributions. 
Because the generalized approximate operators were defined as the generalization of those in conventional particle methods, we could set an optimization problem under wider conditions of parameters than those imposed in conventional particle methods. 
Consequently, the reference weight functions that served as the solution to the optimization problem were different from reference weight functions typically used in conventional particle methods; improvements of accuracy for non-uniform particle distributions were observed through numerical results of the truncation errors and driven cavity flow. 

Finally, we developed the explicit particle method for incompressible Navier--Stokes equations with free surface effects. 
We modified the evaluation of pressure and approximate gradient operator in the explicit particle method to prevent the particle concentrations around the free surface becoming dense. 
We applied the explicit particle method with these modifications to the dam break flow and confirmed a smooth pressure distribution as well as agreement of the time histories of pressure with the experimental results. 

As future work, we will investigate the stability and convergence of the particle methods mathematically. 
Moreover, we will develop particle methods with convergence under more practical conditions such as a that involving fixing the number of neighbor particles $(m=1)$.

\section*{Acknowledgment}
This study was partly supported by priority project 3 for the Post-K Computer entitled ``Sophisticated numerical analysis of diverse earthquake and tsunami disaster scenarios''. 

\appendix

\section{Notation}
\label{sec:appendix_notation}
First, we summarize the computational rules of the multi-index. 
Let $\alpha=(\alpha_1,\alpha_2,\dots,\alpha_\dim)$ be the $\dim$th multi-index. 
For a vector $x\in\dRd$, we denote the $k$th element of $x$ as $\coord{x}{k}$. 
Then, that operations for the multi-index are defined by
\begin{align}
|\alpha|&=\sum_{j=1}^\dim \alpha_j,
\\
x^\alpha&=\prod_{j=1}^\dim (\coord{x}{j})^{\alpha_j},\qquad x\in\dRd,
\\
\alpha! & = \prod_{i=1}^\dim \alpha_i!. 
\end{align}
Let $D^\alpha$ be the differential operator defined by
\begin{align}
D^\alpha=\left(\frac{\partial}{\partial \coord{x}{1}}\right)^{\alpha_1} \left(\frac{\partial}{\partial \coord{x}{2}}\right)^{\alpha_2} \cdots \left(\frac{\partial}{\partial \coord{x}{\dim}}\right)^{\alpha_\dim}, 
\end{align}
where $D^\alpha v=v$ if $|\alpha|=0$. 

Next, we introduce some functional spaces. 
For a set $S\subset\dRd\,(\dim\in\dN)$, let $C(\overline{S})$ be the space of real continuous functions defined in $\overline{S}$, where $\overline{S}$ is the closure of $S$. 
The norm of $C(\overline{S})$ is defined by
\begin{equation}
\brn{\phi}_{C(\overline{S})} := \max_{x\in\overline{S}}\BrA{\phi(x)}. 
\end{equation}
For an open set $S$ and positive integer $k$, let $C^k(\overline{S})$ be the space of functions in $C(\overline{S})$ with derivatives up to the $k$th order. 
The norm of $C^k(\overline{S})$ is defined as
\begin{equation}
\brn{\phi}_{C^k(\overline{S})} := \max_{|\alpha|\leq k}\brn{D^\alpha \phi}_{C(\overline{S})}. 
\end{equation}
Here, $\alpha$ is the multi-index. 
For a functional space $X(\overline{S})$, let $C([0,T];X(\overline{S}))$ be the space of functions on $\overline{S}\times[0,T]$ satisfying 
\begin{equation}
\brn{\phi}_{C([0,T];X(\overline{S}))} := \max_{t\in[0,\tmax]}\brn{\phi(\cdot,t)}_{X(\overline{S})}<+\infty. 
\end{equation}

\section{Description of approximate operators in SPH and MPS using generalized approximate operators}
\label{sec:appendix_Derivation_conventional}
We show that the generalized approximate operators \eqref{approx_intp_01}--\eqref{approx_lap_01} denote approximate operators in SPH and MPS if their parameters are selected appropriately. 
Let $\ws\in\FsWeightFunc$ be a reference weight function such that
\begin{align}
&C_0(\ws)=\int_{\dRd} \ws(|x|)\dx = \int_{\dRd} \wsh(|x|)\dx = 1,
\label{cond:ws_unity}\\
&\WSPHDriv(r)<0,\qquad 0<r<1,
\end{align}
where $\WSPHDriv$ is the first derivative of $\ws$. 
Then, in SPH, the interpolant $\IntpSPH$, approximate gradient operator $\GradAppSPH$, and approximate Laplace operator $\LapAppSPH$ are defined as
\begin{align}
\IntpSPH \f_\i&\deq \sum_{\j=1}^\N\dfrac{m_\j}{\rho_\j}\f_{\j}\wsh(|\Pt{\j}{}-\Pt{\i}{}|),
\label{approx_sph_intp_01}
\\
\GradAppSPH \f_\i&\deq \sum_{\j\neq\i}\dfrac{m_\j}{\rho_\j}(\f_\j-\f_\i)\nabla \wsh (|\Pt{\j}{}-\Pt{\i}{}|),
\label{approx_sph_grad_01}
\\
\LapAppSPH \f_\i &\deq 2 \sum_{\j\neq\i}\dfrac{m_\j}{\rho_\j}\dfrac{\f_\i-\f_\j}{|\Pt{\j}{}-\Pt{\i}{}|}\frac{\Pt{\j}{}-\Pt{\i}{}}{|\Pt{\j}{}-\Pt{\i}{}|} \cdot\nabla \wsh (|\Pt{\j}{}-\Pt{\i}{}|), 
\label{approx_sph_lap_01}
\end{align}
respectively. 
Here, $m_j$ and $\rho_j$ are positive parameters referred to as the particle mass and particle density, respectively. 
The particle volume set $\PvSet$ is given by $\PvSet=\SetNd{\Pv{\i}=m_{\i}/\rho_{\i}}{\i=1,\dots,\N}$. 
Then, from \eqref{cond:ws_unity}, the generalized interpolant \eqref{approx_intp_01} with $\wi=\ws$ is equivalent to the interpolant of SPH \eqref{approx_sph_intp_01}. 
From 
\begin{equation}
-\int_{\dRd} |x| \WSPHDriv(|x|)\dx
= \int_{\dRd} x\cdot\nabla \ws(|x|)\dx = \int_{\dRd} (\nabla\cdot x) \ws(|x|)\dx = \dim\int_{\dRd} \ws(|x|)\dx = \dim, 
\label{eqn:w_SPH_unity}
\end{equation}
the generalized approximate gradient operator \eqref{approx_grad_01} with $\wg=-\WSPHDriv$ is equivalent to the approximate gradient operator of SPH \eqref{approx_sph_grad_01}. 
Moreover, from \eqref{eqn:w_SPH_unity}, the generalized approximate Laplace operator \eqref{approx_lap_01} with
\begin{equation}
\wl(r) = -\dfrac{1}{r}\WSPHDriv(r)
\label{set:w_lap_SPH}
\end{equation}
is equivalent to the approximate Laplace operator of SPH \eqref{approx_sph_lap_01}. 

Let $\wm\in\FsWeightFunc$ be a reference weight function defined by \eqref{def:w_MPS}. 
A weight function $\wmh$ is set by \eqref{def:wh}. 
Then, in MPS, the approximate gradient operator $\GradAppMPS$ and approximate Laplace operator $\LapAppMPS$ are defined as
\begin{align}
\GradAppMPS \f_\i &\deq \frac{\dim}{n_0} \sum_{\j\neq\i}\dfrac{\f_\j-\f_\i}{|\Pt{\j}{}-\Pt{\i}{}|}\frac{\Pt{\j}{}-\Pt{\i}{}}{|\Pt{\j}{}-\Pt{\i}{}|} \wh^{\rm MPS} (|\Pt{\j}{}-\Pt{\i}{}|),
\label{approx_mps_grad_01}\\
\LapAppMPS \f_\i &\deq \frac{2\dim}{n_0\lambda_0} \sum_{\j\neq\i}(\f_\j-\f_\i) \wh^{\rm MPS} (|\Pt{\j}{}-\Pt{\i}{}|), 
\label{approx_mps_lap_01}
\end{align}
respectively. 
Here, $n_0$ and $\lambda_0$ are parameters that depend on both $\wm$ and $\h$. 
In general, $\lambda_0$ is given by $\lambda_0=C_2(\wmh)$. 
Then, the particle volume set $\PvSet$ is given by $\PvSet=\SetNd{\Pv{\i}=C_0(\wmh)/n_0}{\i=1,\dots,\N}$. 
Further, the generalized approximate gradient operator \eqref{approx_grad_01} with
\begin{equation}
\wg(r)=\dfrac{1}{r}\wm(r)
\label{set:w_grad_MPS}
\end{equation}
is equivalent to the approximate gradient operator of MPS \eqref{approx_mps_grad_01}. 
Furthermore, the generalized approximate Laplace operator \eqref{approx_lap_01} with $\wl=\wm$ is equivalent to the approximate Laplace operator of MPS \eqref{approx_mps_lap_01} with $\lambda_0=C_2(\wmh)$. 

\section{Order estimates of approximate pressure}
\label{sec:appendix_proof_discrete_density}
Here, we derive the order estimate \eqref{eq:N-S_Comp:cont_app}. 
We assume $\brnnd{\VelComp}_{C^1([0,\tmax];C^3(\DmH))}<\infty$. 
We arbitrarily set $\k=1,2,\dots,\tstepmax$. 
Let $t\in[\tkn,\tkp]$.  
Then, by the chain rule, we have
\begin{align}
\CompCof^{2}\MDerivComp{\PresCompApp{\k}}(x,t)
&=\dfrac{\Dens}{C_0(\w)}\int_{\DmH}\MDerivComp{\wh}(|\XComp{\k}(y,t)-\XComp{\k}(x,t)|) \dy
\nn\\
&=\dfrac{\Dens}{C_0(\w)}\int_{\DmH}\brm{\VelComp(x,t)-\VelComp(y,t)}\cdot\grad\wh(|\XComp{\k}(y,t)-\XComp{\k}(x,t)|) \dy. 
\label{appendixB:densapp_01}
\end{align}
Further, by Taylor expansion and using our assumptions, we get
\begin{align}
\CompCof^{2}\MDerivComp{\PresCompApp{\k}}(x,t)
&=\dfrac{\Dens}{C_0(\w)}\int_{\DmH}\brm{\VelComp(x,t)-\VelComp(y,t)}\cdot\grad\wh(|y-x|) \dy+\order(\tau\h^{-1})
\nn\\
&=-\dfrac{\Dens}{C_0(\w)}\int_{\DmH}\brm{(y-x)\cdot\nabla}\VelComp(x,t)\cdot\grad\wh(|y-x|)\dy
\nn\\
&\quad-\dfrac{\Dens}{2C_0(\w)}\int_{\DmH}\brm{(y-x)\cdot\nabla}^2\VelComp(x,t)\cdot\grad\wh(|y-x|)\dy+ \order(\tau\h^{-1}+\h^2).
\label{appendixB:densapp_02}
\end{align}
Using the multi-indices $\alpha$ and $\beta$, we have
\begin{multline}
	\int_{\DmH}\brm{(y-x)\cdot\nabla}\VelComp(x,t)\cdot\grad\wh(|y-x|)\dy
	\\
	= -\sum_{|\alpha|=1, |\beta|=1}D^{\alpha}\VelComp(x,t)^{\,\beta} \int_{\DmH}\dfrac{(y-x)^{\alpha+\beta}}{|y-x|}\dfrac{y-x}{|y-x|}\cdot\nabla\wh(|y-x|)\dy. 
\end{multline}
If $\alpha=\beta$, then, by the Gauss--Green theorem and considering 
\begin{equation}
C_0(\wh)=C_0(\w),\qquad \w\in\FsWeightFunc,
\label{appendixB:chara:Ck}
\end{equation}
we have
\begin{align}
\int_{\DmH}\dfrac{(y-x)^{\alpha+\beta}}{|y-x|}\dfrac{y-x}{|y-x|}\cdot\wh(|y-x|)\dy 
&=\int_{\DmH}\dfrac{\brm{(y-x)^{\alpha}}^2}{|y-x|}\dfrac{y-x}{|y-x|}\cdot\wh(|y-x|)\dy
\nn\\
&=-\dfrac{1}{\dim}\int_{\DmH}(y-x)\cdot\nabla\wh(|y-x|)\dy
\nn\\
&=\int_{\DmH}\wh(|y-x|)\dy
\nn\\
&=C_0(\w). 
\end{align} 
If $\alpha\neq\beta$, then, by the symmetry of the integrated function, we have
\begin{equation}
\int_{\DmH}\dfrac{(y-x)^{\alpha+\beta}}{|y-x|}\dfrac{y-x}{|y-x|}\cdot\wh(|y-x|)\dy = 0. 
\end{equation}
Thus, we obtain
\begin{equation}
\int_{\DmH}\brm{(y-x)\cdot\nabla}\VelComp(x,t)\cdot\grad\wh(|y-x|)\dy = C_0(\w)\diver\VelComp(x,t). 
\label{appendixB:densapp_03}
\end{equation}
Moreover, by the symmetry of the integrated function, the second term on right-hand side in \eqref{appendixB:densapp_02} becomes
\begin{multline}
\dfrac{1}{2}\int_{\DmH}\brm{(y-x)\cdot\nabla}^2\VelComp(x,t)\cdot\grad\wh(|y-x|)\dy \\
=\dfrac{1}{2}\sum_{|\alpha|=2, |\beta|=1}\brm{D^{\alpha}\VelComp(x,t)}^{\beta} \int_{\DmH}\dfrac{(y-x)^{\alpha+\beta}}{|y-x|}\deriv{1}{r}{}\wh(|y-x|)\dy
=0. 
\label{appendixB:densapp_04}
\end{multline}
Therefore, by \eqref{appendixB:densapp_01}, \eqref{appendixB:densapp_02}, \eqref{appendixB:densapp_03}, and \eqref{appendixB:densapp_04}, we obtain
\begin{equation}
\CompCof^2\MDerivComp{\PresCompApp{\k}}(x,t)+\Dens\diver\VelComp(x,t) =\order(\tau\h^{-1}+\h^2),\qquad x\in\Dm,\quad t\in[\tkn,\tkp].
\end{equation}
	
 \newcommand{\noop}[1]{}

\end{document}

%% file: new_command.tex
\def\deq{\mathrel{\mathop:}=}
\newcommand{\midd}{\;\middle|\;}
\renewcommand{\dim}{d}
\newcommand{\dN}{\mathbb{N}}

\newcommand{\dR}{\mathbb{R}}
\newcommand{\dRd}{\dR^\dim}

\newcommand{\dZ}{\mathbb{Z}}
\newcommand{\gausssymbol}[1]{\lfloor #1 \rfloor}
\newcommand{\ra}{\rightarrow}

\newcommand{\nn}{\nonumber}
\newcommand{\brs}[1]{\left(#1\right)}
\newcommand{\brm}[1]{\left\{#1\right\}}

\newcommand{\BrA}[1]{\left|#1\right|}
\newcommand{\brn}[1]{\left\|#1\right\|}
\newcommand{\brnnd}[1]{\|#1\|}
\newcommand{\ds}{\displaystyle}
\newcommand{\dx}{{\rm \,d}x}
\newcommand{\dy}{{\rm \,d}y}
\newcommand{\coord}[2]{{#1}^{(#2)}}
\newcommand{\meas}[1]{\BrA{#1}} 
\newcommand{\Set}[2]{\left\{#1\midd#2\right\}}
\newcommand{\SetNd}[2]{\{#1\mid#2\}}

\newcommand{\Dm}{\Omega}
\newcommand{\bd}{\Gamma}
\newcommand{\Vel}{u}
\newcommand{\Pres}{p}
\newcommand{\force}{f}

\newcommand{\Dens}{\rho}
\newcommand{\Velini}{\Vel_0}
\newcommand{\Velbd}{\Vel_{\bd}}
\newcommand{\tmax}{T}
\newcommand{\grad}{\nabla}
\newcommand{\lap}{\Delta}
\newcommand{\diver}{\nabla\cdot}

\renewcommand{\H}{H}
\newcommand{\timeindex}[1]{#1}
\newcommand{\DmH}{\Dm_{\H}}
\renewcommand{\i}{i}
\renewcommand{\j}{j}
\newcommand{\PtChar}{x}
\newcommand{\Pt}[2]{\PtChar_{#1}^{\timeindex{#2}}}
\newcommand{\N}{N}
\newcommand{\PtSetChar}{\mathcal{X}}
\newcommand{\PtSet}{\PtSetChar_\N}
\newcommand{\PvCher}{\omega}
\newcommand{\Pv}[1]{\PvCher_{#1}}
\newcommand{\PvSet}{\mathcal{V}_\N}
\newcommand{\w}{w}
\newcommand{\h}{h}
\newcommand{\wh}{\w_\h}
\newcommand{\f}{\phi}
\newcommand{\Intp}[1]{\Pi_{\h}^{\timeindex{#1}}}
\newcommand{\GradApp}[1]{\grad_{\h}^{\timeindex{#1}}}

\newcommand{\LapApp}[1]{\lap_{\h}^{\timeindex{#1}}}

\newcommand{\IntpSPH}{\Pi_{\h}^{\rm SPH}}
\newcommand{\GradAppSPH}{\grad_{\h}^{\rm SPH}}
\newcommand{\LapAppSPH}{\lap_{\h}^{\rm SPH}}
\newcommand{\GradAppMPS}{\grad_{\h}^{\rm MPS}}
\newcommand{\LapAppMPS}{\lap_{\h}^{\rm MPS}}
\newcommand{\IntpCof}{C_{\Pi}}
\newcommand{\GradAppCof}{C_{\grad}}
\newcommand{\LapAppCof}{C_{\lap}}
\newcommand{\FsWeightFunc}{\mathcal{W}}

\newcommand{\wi}{\w^{\Pi}}
\newcommand{\wg}{\w^{\grad}}
\newcommand{\wl}{\w^{\lap}}
\newcommand{\wih}{\wi_\h}
\newcommand{\wgh}{\wg_\h}
\newcommand{\wlh}{\wl_\h}

\newcommand{\bdh}{\bd_{\H}}
\newcommand{\Veliniex}{\Vel_{0\!\!\;,\,\H}}
\newcommand{\Velbdex}{\Vel_{\bd\!\!\;,\,\H}}
\newcommand{\Dt}{\tau}
\newcommand{\Dtmax}{\Dt_{\max}}
\newcommand{\tstep}{k}
\renewcommand{\k}{\tstep}
\newcommand{\tstepmax}{K}
\newcommand{\timecher}{t}

\newcommand{\tk}{\timecher^{\timeindex{\tstep}}}
\newcommand{\tkn}{\timecher^{\timeindex{\tstep-1}}}
\newcommand{\tkp}{\timecher^{\timeindex{\tstep+1}}}
\newcommand{\ptt}[2]{\PtChar_{#1}^{\timeindex{#2}}}
\newcommand{\PtP}[2]{\PtChar_{#1}^{\ast\!\!\;,\,\timeindex{#2}}}
\newcommand{\PtSetTime}[1]{\PtSetChar_\N^{\timeindex{#1}}}
\newcommand{\PtSetTent}[1]{\PtSetChar_\N^{\ast\!\!\;,\,\timeindex{#1}}}
\DeclareMathOperator{\esssup}{ess\,\sup}
\newcommand{\VelAppChar}{\mathrm{\Vel}}
\newcommand{\VelApp}[2]{\VelAppChar_{#1}^{\timeindex{#2}}}
\newcommand{\VelAppP}[2]{\VelAppChar_{#1}^{\ast\!\!\;,\,\timeindex{#2}}}
\newcommand{\PresAppChar}{\mathrm{\Pres}}
\newcommand{\PresApp}[2]{\PresAppChar_{#1}^{\timeindex{#2}}}
\newcommand{\PresAppP}[2]{\PresAppChar_{#1}^{\ast\!\!\;,\,\timeindex{#2}}}

\newcommand{\IndexSet}[2]{\Lambda^{\timeindex{#2}}(#1)}

\newcommand{\mderivD}{{\rm D}}
\newcommand{\mderiv}[3]{
	\ifnum #1=1	\frac{\mderivD#3}{\mderivD{#2}}%
	\else{\frac{\mderivD^{#1}#3}{\mderivD^{#1}{#2}} %
	}\fi%
}

\newcommand{\IntpMod}[1]{\widehat{\Pi}_{\h}^{\timeindex{#1}}}
\newcommand{\GradAppPlus}[1]{\grad_{\h,+}^{\timeindex{#1}}}
\newcommand{\GradAppPlustent}[1]{\grad_{\h,+}^{\ast\!\!\;,\,\timeindex{#1}}}

\newcommand{\CompCof}{\varepsilon}
\newcommand{\VelComp}{\Vel_{\CompCof}}
\newcommand{\PresComp}{\Pres_{\CompCof}}

\newcommand{\Xcher}{X}
\newcommand{\XComp}[1]{\Xcher^{#1}_{\CompCof}}
\newcommand{\Ddel}{{\rm D}_{t, \CompCof}}
\newcommand{\MDerivComp}[1]{\Ddel\,#1}

\newcommand{\pderivD}{{\rm \partial}}
\newcommand{\pderiv}[3]{
	\ifnum #1=1	\frac{\pderivD#3}{\pderivD{#2}}%
	\else{\frac{\pderivD^{#1}#3}{\pderivD^{#1}{#2}} %
	}\fi%
}

\newcommand{\VelCompDt}[1]{\Vel_{\CompCof, \Dt}^{\timeindex{#1}}}
\newcommand{\PresCompDt}[1]{\Pres_{\CompCof, \Dt}^{\timeindex{#1}}}

\newcommand{\VelCompDtP}[1]{\Vel_{\CompCof, \Dt}^{\ast\!\!\;,\,\timeindex{#1}}}

\newcommand{\XCompDt}[1]{\Xcher^{\timeindex{#1}}_{\CompCof, \Dt}}
\newcommand{\XCompDtP}[1]{\Xcher^{\ast\!\!\;,\,\timeindex{#1}}_{\CompCof, \Dt}}

\newcommand{\order}{\mathcal{O}}

\newcommand{\PresCompApp}[1]{\widetilde{\Pres}^{#1}_{\CompCof}}

\newcommand{\SPHCher}{\rm SPH}
\newcommand{\ws}{\w^{\SPHCher}}
\newcommand{\wsh}{\ws_{\h}}
\newcommand{\WSPHDriv}{\dot{\w}^{\SPHCher}}
\newcommand{\wm}{\w^{\rm MPS}}
\newcommand{\wmh}{\wm_{\h}}

\newcommand{\DmHo}{\overline{\Dm}_{H}}

\newcommand{\Pertab}{\epsilon}
\newcommand{\PertabMax}{\Pertab_{\max}}

\newcommand{\Dx}{\Delta x}
\newcommand{\PresAppave}{\overline{\PresAppChar}}
\newcommand{\PresAppk}{\PresAppChar^{\timeindex{\tstep}}}

\newcommand{\PresAppavek}{\PresAppave^{\timeindex{\tstep}}}
\newcommand{\derivD}{{\rm d}}
\newcommand{\deriv}[3]{
	\ifnum #1=1	\frac{\derivD#3}{\derivD{#2}}%
	\else{\frac{\derivD^{#1}#3}{\derivD^{#1}{#2}} %
	}\fi%
}
